\documentclass{amsart}
\usepackage[utf8]{inputenc}

\usepackage{tikz}

\newtheorem{theorem}{Theorem}
\newtheorem{lemma}[theorem]{Lemma}
\newtheorem*{lemma*}{Lemma}

\newtheorem{claim}[theorem]{Claim}
\theoremstyle{definition}
\newtheorem{definition}[theorem]{Definition}
\newtheorem{remark}[theorem]{Remark}
\newtheorem{notation}[theorem]{Notation}
\newtheorem{example}[theorem]{Example}
\newtheorem{fact}[theorem]{Fact}
\newtheorem{question}{Question}

\newcommand{\R}{\mathbb{R}}
\newcommand{\sm}{\setminus}

\newcommand{\cH}{\mathcal{H}}

\newcommand{\interior}[1]{{#1}^\circ}
\newcommand{\closure}[1]{\overline{#1}}
\newcommand{\concat}{{}^\smallfrown}

\usepackage{amssymb}

\newcommand{\I}{{\boldsymbol{I}}}
\newcommand{\II}{{\boldsymbol{I}\kern-0.05cm\boldsymbol{I}}}
\newcommand{\res}{\mathord \upharpoonright}

\newcommand{\abs}[1]{\left | #1 \right |}
\newcommand{\setof}[1]{\left\{#1\right\}}
\newcommand{\suchthat}{~\colon~}
\newcommand{\sgn}{\text{sgn}}
\newcommand{\ac}{\mathsf{AC}}
\newcommand{\dc}{\mathsf{DC}}
\newcommand{\Q}{\mathbb{Q}}
\newcommand{\ad}{\mathsf{AD}}
\newcommand{\zf}{\mathsf{ZF}}
\newcommand{\zfc}{\mathsf{ZFC}}

\newcommand{\bp}[2]{[#1,#2]_P}
\DeclareMathOperator{\lspan}{span}

\newcommand{\ignore}[1]{}
\title{The no-$\beta$ McMullen Game and the perfect set property}

\newcommand{\authorlogan}{
\author{Logan Crone}
\address{Logan Crone, University of North Texas, Department of Mathematics, 
1155 Union Circle \#311430, Denton, TX 76203-5017, USA}
\email{logancrone@my.unt.edu}
}

\newcommand{\authorlior}{
\author{Lior Fishman}
\address{University of North Texas, Department of Mathematics, 1155 Union Circle \#311430, Denton, TX 76203-5017, USA}
\email{lior.fishman@unt.edu}
}

\newcommand{\authorsteve}{
\author{Stephen Jackson}
\address{Stephen Jackson, University of North Texas, Department of Mathematics, 1155 Union Circle \#311430, Denton, TX 76203-5017, USA}
\email{jackson@unt.edu}
}
\thanks{The third author was partially supported by NSF grant DMS-1800323.}

\begin{document}

\authorlogan
\authorlior
\authorsteve



\maketitle

\begin{abstract}
Given a target set $A\subseteq  \R^d$ and a real number $\beta\in (0,1)$, McMullen \cite{McMullen2010}
introduced the notion of $A$ being an absolutely $\beta$-winning set.
This involves a two player game which we call the $\beta$-McMullen game. We consider the version of this game
in which the parameter $\beta$ is removed, which we call the no-$\beta$ McMullen game.
More generally, we consider the game with respect to arbitrary norms on $\R^d$, and even more generally
with respect to general convex sets. We show that for strictly convex sets in  $\R^d$, polytopes in $\R^d$,
and general convex sets in $\R^2$, that player $\I$ wins the no-$\beta$ McMullen game iff
$A$ contains a perfect set and player $\II$ wins iff $A$ is countable. 
So, the no-$\beta$ McMullen game is equivalent to the
perfect set game for $A$ in these cases. The proofs of these results use a  connection 
between the geometry of the game and techniques from logic. 
Because of the geometry of this game, this result has strong implications for the geometry
of uncountable sets in $\R^d$. We also present an example of a compact, convex set in $\R^3$ to which our methods do not
apply, and also an example due to D.\ Simmons of a closed, convex set in $\ell_2(\R)$
which illustrate the obstacles in extending the results further. 
\end{abstract}

\section{Introduction} \label{sec:intro}

In McMullen's original game we have a target set $A\subseteq \R^d$ and a parameter $\beta \in (0,1/3)$.
The two players, which we call $\I$ and $\II$, alternate playing closed balls in $\R^d$, with $\I$ making moves
$B_{2n}=B(x_{2n},\rho_{2n})$ and $\II$ making moves $B_{2n+1}=B(x_{2n+1},\rho_{2n+1})$. The rules of the game are that
$B_{2n+2} \subseteq B_{2n}\sm B_{2n+1}$ (this is only a requirement on $\I$, note that we do not require that
$B_{2n+1} \subseteq B_{2n}$). Also, we require that $\rho_{2n+2} \geq \beta \rho_{2n}$ and
$\rho_{2n+1} \leq \beta \rho_{2n}$. In a run of the game following the rules, a unique point
$x= \lim_n x_n$ is determined (so $\{x\}= \bigcap_n B_{2n}$). Then $\I$ wins the run of the game
iff $x \in A$.

In the no-$\beta$ McMullen game, we change only the rule restricting the radii.
We instead require only that $\rho_{2n}>0$ and that $\rho_{2n+1} < \rho_{2n}$.  In this
variation of the original game, we may not have that $\rho_n \to 0$,
and so we declare $\I$ the winner just in case $\bigcap_{n} B_{2n} \cap A \neq \emptyset$.

More generally, given a convex, compact set $P \subseteq \mathbb{R}^d$,
we define another game in which the players alternate playing translated and
scaled copies of $P$ (say by playing pairs $(x_n, \rho_n)$ to play the set $P_n = x_n + \rho_n P$).
Again we require that $P_{2n+2} \subseteq P_{2n} \setminus P_{2n+1}$ and $\rho_{2n} > 0$,
and $\rho_{2n+1} < \rho_{2n}$.
If the players follow these rules, we declare $\I$ the winner if $\bigcap_{n} P_{2n} \cap A \neq \emptyset$.
A particular instance of this would be when the convex compact sets $P$ are the closed
unit balls with respect to some norm $\| \cdot \|$ on $\R^n$. This gives the no-$\beta$
McMullen game as a special instance of this more general game.
We let $G_{P}(A)$ denote this no-$\beta$ McMullen game just defined. 
See Figure~\ref{fig:nb} for an illustration of the moves of the game.

We work throughout in the base theory $\zf +\dc$. This will allow our main results to be 
presented both in the $\zfc$ context for Borel target sets and in the $\ad$ context 
for arbitrary target sets. In fact, our main Theorems~\ref{thma}, \ref{thmb}, \ref{thmc}, \ref{geomthm}
hold as stated for arbitrary sets with no determinacy hypothesis but, for example, to apply Theorem~\ref{geomthm}
to general uncountable sets in $\R^d$ we need have the perfect set property, a consequence 
of $\ad$.

\newcommand{\polygon}[4][draw=black]{\draw [#1](#2*-1+#3, #2*0+#4) --
  (#2*0+#3, #2*2+#4) -- (#2*1+#3, #2*1+#4) -- (#2*1.5+#3, #2*-1+#4) --
  (#2*-0.5+#3, #2*-1+#4) -- (#2*-1+#3, #2*0+#4);}

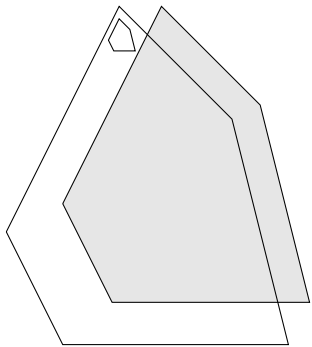
\begin{figure}[h]
  \begin{tikzpicture}[scale=0.75]
    \polygon[draw=black, fill=gray!20!white]{1.75}{0.75}{0.5}
    \polygon[draw=black, fill=none]{2}{0}{0}
    \polygon[draw=black, fill=none]{0.19}{0}{3.4}
  \end{tikzpicture}
  \caption{The no-$\beta$ McMullen game} 
\label{fig:nb}
\end{figure}

A fundamental theorem of Martin \cite{Martin1975} (see also \cite{Martin1985}) asserts that
every Borel game on an arbitrary set $X$ is determined (without assuming $\ac$ one only
gets a {\em quasistrategy} for one of the players; see \S 6F of \cite{MoschovakisBook} for further discussion).
It follows easily that for any Borel target set $A \subseteq X$, where $X$ is a separable Banach space,
that the no-$\beta$ game is determined (or without $\ac$, quasidetermined). 
In Remark~\ref{acremark} we note that for the no-$\beta$ McMullen game that in fact we do not need
full $\ac$ but just $\dc$ to get a strategy in the case that the target set $A$ is Borel.
In Remark~\ref{detremark} we will note that for any $P \subseteq \R^d$ for which we can prove our main theorem
we have that the no-$\beta$ game $G_P(A)$ is determined for all $A\subseteq \R^d$ from just the
axiom of determinacy, $\ad$, even though it is a real game.

\begin{remark} \label{acremark}
For Borel target sets $A\subseteq X$ ($X$ any separable Banach space), the determinacy
of the game $G_P(A)$ follows from Martin's theorem if we assume $\ac$. Assuming only
$\dc$, the determinacy of integer games, we also get the determinacy of $G_P(A)$ as follows. 
First, it is not hard to see that the game $G_P(A)$ is equivalent to the version $G'_P(A)$
in which $\II$ must play scaled copies $s+tP$ with $s \in \Q^d$, $t\in \Q$ (``equivalent'' here means that
if one of the players has a winning strategy in one of the games, then that same player has a winning
strategy in the other game). Without $\ac$, Martin's theorem gives that one of the players has a
winning  quasistrategy  for $G'_P(A)$. If $\II$ has a winning quasistrategy for $G'_P(A)$,
then (since $\II$'s moves are coming from a countable set) $\II$ actually has a winning strategy
for $G'_P(A)$ and thus a winning strategy for $G_P(A)$ (using a fixed 
wellordering of the countable set of $\II$'s moves,
$\II$ plays at each round the least move consistent with the quasistrategy).
If $\I$ has a winning quasistrategy in $G'_P(A)$ then, using $\dc$ and the fact that $\II$'s moves are coming
from a countable set, we can turn the quasistrategy into a strategy. Similarly, 
if $A$ is closed, then the game $G_P(A)$ is easily a closed game and the by Gale-Stewart and so by the 
argument above the game $G_P(A)$ is determined. 
\end{remark}

Removing the $\beta$ restrictions from the McMullen game gives the game some of the flavor
of the Banach-Mazur game. In \cite{FSR2016} an intermediate form of the game was considered where
the $\beta$ requirement on player $\I$ was removed, but the $\beta$ requirement for player $\II$
is kept. In this case it was shown that the game is equivalent to the perfect set game.
Recall the perfect set game, which was introduced in \cite{Davis1964} (see also \S 6A of \cite{MoschovakisBook}
for further discussion), is a game defined for any set $A\subseteq X$, for $X$ a Polish space, and has the property that
$\II$ wins iff $A$ is countable, and $\I$ wins iff $A$ contains a perfect set. Removing the $\beta$
requirements from player $\II$ as well superficially seems to give player $\II$ much more control.
It is perhaps somewhat surprising then that our main results say that the no-$\beta$ game
is still equivalent to the perfect set game, at least for the classes of $P$
for which we have been able to prove the theorem. This result has consequences for the geometric structure
of sets in $\R^d$, which we state below.

We say a norm $\| \cdot \|$ on $\R^d$ with closed unit ball $P$ is {\em strictly convex} (or we say the
convex compact set $P$ is strictly convex) if for every $x,y \in P$ every point of the line
segment $\ell$ between $x$ and $y$ except for $x$ and $y$ lies in the interior of $P$.
Equivalently, $P$ is strictly convex if there is no proper (i.e., not a point) line
segment contained in $\partial P$.

In the following Theorems, when we say that a game $G$ is ``equivalent to the perfect set game''
we mean simply that $\I$ has a winning strategy for the game $G$  iff the target set $A$ contains a perfect set, and $\II$
has a winning strategy for $G$ iff $A$ is countable.

We note that the hypotheses of the following Theorems~\ref{thma}, \ref{thmb}, \ref{thmc} 
all imply that the convex set $P\subseteq \R^d$ has a non-empty interior. This is necessary as the 
theorem does not hold in general without this assumption.

The first class of $P$ for which we are able to prove our result is the class of
strictly convex sets.

\begin{theorem} \label{thma}
Let  $P$ be a non-degenerate (i.e., more than one point)
compact, strictly convex set in $\R^d$.  Then the no-$\beta$
McMullen game $G_P(A)$ with respect to $P$ and target set $A$ is equivalent to the perfect set game.
\end{theorem}

In particular, if $\| \cdot \|$ is a strictly convex norm on $\R^n$ then 
the no-$\beta$ McMullen game with respect to the closed unit ball for the norm
is equivalent to the perfect set game.

The second class of $P \subseteq \R^d$ for which we can establish our result
are the convex $d$-polytopes (the compact convex sets  $P$ which are the intersections
of finitely many closed half-spaces and have non-empty interior). For example, convex polygons in $\R^2$ and
convex polyhedra in $\R^3$.

\begin{theorem} \label{thmb}
Let $P$ be a compact, convex $d$-polytope in $\R^d$. Then the no-$\beta$ McMullen game $G_P(A)$
is equivalent to the perfect set game.
\end{theorem}

Finally, we are able to establish the result for arbitrary compact, convex set in $\R^2$.

\begin{theorem} \label{thmc}
Let $P$ be a compact, convex set in $\R^2$ with non-empty interior.
Then the no-$\beta$ McMullen game $G_P(A)$ is equivalent to the perfect set game.
\end{theorem}

In showing the equivalence of the game $G_P(A)$ with the perfect set game, two 
of the four implications are easy. If $A$ is countable, then easily $\II$ 
as a winning strategy by consecutively deleting the members of $A$ 
in each round. If $\I$ has a winning strategy, then $A$ contains a perfect set. 
In fact, this holds in a general Banach space which we state in the next fact.

\begin{fact} \label{psf}
Let $P$ be a bounded, closed, convex set in a Banach space $X$, and let $A\subseteq X$. 
If $\I$ has a winning strategy for the no-$\beta$ McMullen game $G_P(A)$, then $A$
contains a perfect set.
\end{fact}

\begin{proof}
Let $\sigma$ be a winning strategy for $\I$ in $G_P(A)$. Without loss of generality
we may assume that $\vec 0 \in P$ (where $\vec 0$ is the zero vector of the Banach space). 
We make use of the following observation. If $Q=y+\delta P \subseteq P\sm (1-\epsilon)P$,
then $\delta \leq \epsilon$. Let $L$ be the line in $X$ through $\vec 0$ and $y$. 
Let $m=|L \cap P|$. Then $|L\cap (1-\epsilon)P|= (1-\epsilon)m$ and 
$|L\cap Q|= \delta m$. The second equality follows as $Q\cap L= Q\cap (y+L)
=(y+\delta P)\cap (y+L)= y+(\delta P\cap L)$. But then 
$m=|L\cap P| \geq |L \cap Q| +|L\cap (1-\epsilon)P| =\delta m+ (1-\epsilon)m$,
and so $\delta \leq \epsilon$.

We build a perfect set contained as $A$ as follows.
Let $\I$ play $P_\emptyset=x_\emptyset+\rho_\emptyset P$ according to $\sigma$. In general suppose we have
$P_s= x_s+ \rho_s P$ defined for some $s \in 2^{<\omega}$. Let $P'_{s \concat 0} \subseteq P_s$
have a scaling factor $> \frac{1}{2}\rho_s$. Let $P_{s\concat 0}= x_{s \concat 0} +\rho_{s \concat 0} P$ be $\I$'s response
to $P'_{s \concat 0}$. By the above observation we have that $\rho_{s \concat 0} < \frac{1}{2}\rho_s$. 
Let $P'_{s \concat 1}= P_{s\concat 0}$, and let $P''_{s \concat 1}$ be $\I$'s response to 
$P'_{s \concat 1}$. Let $P'''_{s\concat 1}\subseteq P''_{s \concat 1}$ have a scaling 
factor greater that one-half of that for $P''_{s \concat 1}$. Finally, let $P_{s \concat 1}$
be $\I$'s response to $P'''_{s\concat 1}$. Again by the observation, $P_{s \concat 1}$
has scaling factor $< \frac{1}{2} \rho_s$. 
So, $P_{s\concat 0}, P_{s \concat 1} \subseteq P_s$,
$P_{s\concat 0}\cap P_{s \concat 1} =\emptyset$, and $\rho_{s \concat 0}$, $\rho_{s \concat 1} < \frac{1}{2}
\rho_s$.
For each $\zeta \in 2^\omega$, let $\{ x_\zeta\} =\bigcap_n P_{\zeta \res n}$. Then $\{ x_\zeta \colon
\zeta \in 2^\omega\}$ is a perfect set contained in $A$. 
 \end{proof}

In the proofs of Theorems~\ref{thma}, \ref{thmb}, and \ref{thmc}, we will actually
show that if $\II$ has a winning strategy in $G_P(A)$ then $A$ is countable.
The following remark shows that this is enough to obtain the full statements of
these theorems.

\begin{remark} \label{detremark}
Let $P\subseteq \R^d$ be as in Theorems~\ref{thma}, \ref{thmb}, or \ref{thmc}.
Suppose we have shown that for any target set $A$ if $\II$ wins $G_P(A)$ then $A$
is countable. If $A$ is countable, easily $\II$ wins $G_P(A)$ and so 
$\II$ wins $G_P(A)$ iff $A$ is countable. If $\I$ has a winning strategy 
then by Fact~\ref{psf}, $A$ contains a perfect set. 
On the other hand, if $A$ contains a perfect set $F$, then by Remark~\ref{acremark} the game $G_P(F)$
is determined. Since $F$ is uncountable, $\II$ cannot have a winning strategy
in $G_P(F)$ by assumption, and so $\I$ has a winning strategy for $G_P(F)$. 
A winning strategy for $G_P(F)$ is easily a winning strategy for $G_P(A)$ as well.

Also, for any target set $A \subseteq \R^d$ 
the determinacy of the game $G_P(A)$ follows
just from $\ad$, even though these are real games. To see this, let $A\subseteq \R^d$.
By $\ad$, either $A$ is countable or else $A$ contains a perfect set. If $A$ is countable,
then as above easily $\II$ wins $G_P(A)$. If $A$ contains a perfect set $F$, 
then also as above $\I$ wins $G_P(F)$ and thus $G_P(A)$. 
\end{remark}

To state the geometric consequences of these results, we introduce the notion
of the derivative of a set $A\subseteq \R^d$ with respect to the set $P$. The notion is
somewhat analogous to the usual Cantor-Bendixson derivative for sets, which is obtained by iterating the operation
of taking limit points (and taking intersections at limit stages).

\begin{definition} \label{derivative}
Let $A\subseteq \R^d$, and $P\subseteq \R^d$ a compact, convex set. Then the derivative
$A'_P$ of $A$ with respect to $P$ is defined by: $A\sm A'_P$ is the union of all
open balls which do not contain a {\em good copy} of $P$. We say $Q= x+ \rho P$ ($Q$
is a translated, scaled copy of $P$) is a good copy of $P$ if
the set of points $B\subseteq A\cap \partial Q$ which are limit points of either $A \cap \interior Q$
or $A \cap \partial Q \cap H$ where $H$ is a supporting hyperplane for $Q$ in $\R^d$,
cannot be covered by a strictly smaller copy of $P$ (i.e., some $x' + \rho' P$ where
$\rho' <\rho$).
\end{definition}

Note that the derivative operation is monotone, that is, if $A\subseteq B$ then
$A'_P \subseteq B'_P$.

As usual, we iterate this notion of derivative, taking intersections at limit ordinals.
This produces a sequence of sets $A=A^0_P \supseteq A^1_P \supseteq \cdots \supseteq A^\alpha_P
\supseteq \cdots \supseteq A^\infty_P$ where $A^\infty_P=(A^\infty_P)'$. 
The usual argument shows that starting from
any set $A \subseteq \R^d$, the derivative stabilizes at a countable ordinal, that is,
there is a countable ordinal $\alpha$ so that $A^\infty_P=A^\alpha_P$. 
As usual, if we start with a closed set $A$, then all of the $A^\alpha_P$ are closed.

Note that for any set $A$, $A'_P \subseteq A'$, where $A'$ here denotes the usual
Cantor-Bendixson derivative. So if $F$ is a closed set and $F^\infty_P=F$
then not only is $F$ perfect in the usual sense, but $F$ has strong geometric regularity
with respect to the shape $P$.

In section~\ref{sec:oa} we get the following geometric 
consequence of the above theorems and a derivative analysis.

\begin{theorem} \label{geomthm}
Suppose $P\subseteq \R^d$ is a compact, convex set such that
the hypothesis of Theorems~\ref{thma}, \ref{thmb}, or \ref{thmc} is satisfied. 
Then if $A \subseteq \R^d$ contains a perfect set then $A^\infty_P \neq \emptyset$. 
\end{theorem}

As a consequence, if $A\subseteq \R^d$ is any uncountable Borel set, then
$A$ contains a perfect set $F$ such that $F'_P=F$. So, for every point $x \in F$ and every
neighborhood $U$ of $X$, $U$ contains a good copy of $P$.  Similarly, if we assume 
$\ad$ then this result holds for any $A\subseteq \R^d$.

\begin{figure}[h] 
\centering
\begin{tikzpicture}[scale=1.5]
  \draw [fill=white, draw=black] (0, 0) circle (1.5cm);
  \draw [fill=gray!20!white] (1.1, 1) -- (-1, 1.1) -- (0.05, -1.49) -- (1.1, 1);
  \draw [fill=white, draw=black] (0,0) circle(0.035cm);
  \draw [fill=black, draw=black] (1.1, 1) circle (0.035cm);
  \draw [fill=black, draw=black] (0.925*1.1, 0.925*1) circle (0.025cm);
  \draw [fill=black, draw=black] (0.85*1.1, 0.85*1) circle (0.025cm);
  \draw [fill=black, draw=black] (0.7*1.1, 0.7*1) circle (0.025cm);
  \draw [fill=black, draw=black] (-1, 1.1) circle (0.035cm);
  \draw [fill=black, draw=black] (0.925*-1, 0.925*1.1) circle (0.025cm);
  \draw [fill=black, draw=black] (0.85*-1, 0.85*1.1) circle (0.025cm);
  \draw [fill=black, draw=black] (0.7*-1, 0.7*1.1) circle (0.025cm);
  \draw [fill=black, draw=black] (0.05, -1.49) circle (0.035cm);
  \draw [fill=black, draw=black] (0.925*0.05, 0.925*-1.49) circle (0.025cm);
  \draw [fill=black, draw=black] (0.85*0.05, 0.85*-1.49) circle (0.025cm);
  \draw [fill=black, draw=black] (0.7*0.05, 0.7*-1.49) circle (0.025cm);
\end{tikzpicture}
\caption{A good copy of the unit ball in $\R^2$} \label{fig:gc}
\end{figure}
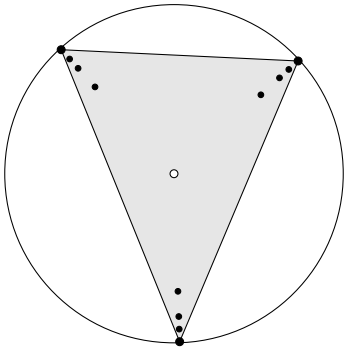

To illustrate this statement with an example, consider the case where $P$ is
the standard closed ball in $\R^2$ (i.e., we use the standard Euclidean norm on $\R^2$).
Let $A \subseteq \R^2$ be any uncountable Borel set. Then $A$ contains a perfect set
$F$ such that for any $x \in F$ and any neighborhood $U$ of $x$,  there is a
good closed ball $B=B(y,r)\subseteq U$ for $F$,
that is, the points of $F\cap \partial B$ which are limit
points of $F\cap \interior{B}$
are not contained in any strict half of $B$. More precisely, these points do not lie strictly
on one side of a line through $y$, or equivalently, $y$ is in the
convex hull of these points. See Figure~\ref{fig:gc}.

We introduce some notation that we will use for the rest of the paper. 
Let $P$
be a compact, convex set in $\R^d$ with non-empty interior.

\begin{definition} \label{def:maximal}
Let $\ell \subseteq P$ be a line segment in $\R^d$. We say $\ell$ is {\em maximal}
with respect to $P$ if there is no translated enlarged copy $\ell'=x+ s\ell$, with $s>1$,
and $\ell'\subseteq P$.
\end{definition}

Note that $\ell$ being maximal is equivalent to saying that $\ell \nsubseteq P'$
for any translated smaller copy $P'=x+sP$, where $s<1$, of $P$.

\begin{remark} \label{rem:lb}
Note the simple but important fact that there is a lower bound $\epsilon(P)>0$
on the lengths of the maximal line segments with respect to $P$. This is immediate from the fact
that $\interior{P}\neq \emptyset$.
\end{remark}

Given distinct points $x, y \in \R^d$, let $\bp x y= w+ sP$ be a scaled copy of $P$
containing both $x, y$ and such that no copy $w'+s'P$ with $s'<s$ contains both
$x$ and $y$. In general, $\bp{x}{y}$ is not unique, and we just fix a choice for it
(this can be done without any form of $\ac$; the relation $R((x,y),(w,s))$
satisfying this definition is a Borel relation with compact sections).  
In the case where $P$ is strictly convex, there is a unique choice for $\bp{x}{y}$. 
In section~\ref{sec:poly}, however,  we will be more particular in our choice of $\bp{x}{y}$.
Note that $x,y$ are always in the boundary of $\bp{x}{y}$, and 
$\ell_{x,y}$ is maximal in $\bp{x}{y}$.

\section{The strictly convex case} \label{sec:sc}

In this section we prove our main result, Theorem~\ref{thma}, in the case where $P$ is a compact,
non-degenerate (at least two points)
strictly convex set in $\R^d$. 
The proof will make use of an idea from logic involving elementary substructures,
and this idea will be used in the other cases as well.

Throughout the rest of this section $P$ will denote a compact, non-degenerate, strictly
convex set in $\R^d$. This implies $\interior{P} \neq \emptyset$.

For $x \in \partial P$, let $\cH_P(x)$ denote the collection of supporting hyperplanes
for $P$ at $x$. That is, $H\in \cH_P(x)$ iff $H \subseteq \R^d$ is an affine hyperplane containing $x$
with $P$ entirely contained (not-strictly) on one side of $H$. If $H$ is an affine hyperplane in $\R^d$
and $u$ a vector in $\R^d$, then by $|H\cdot u|$ we mean the absolute value of the dot-product 
$|u\cdot n_H|$ where $n_H$ is a unit normal vector to $H$. If $\ell$ is a line segment,
we also use $| \frac{\ell}{\| \ell\| } \cdot H |$ to denote the dot-product where 
we interpret $\ell$ as a vector (the sign ambiguity doesn't matter due to the absolute value).

\begin{definition}
Let $x \in \partial P$. Let $M_P(x)$ be the collection of line segments $\ell$ which are maximal with
respect to $P$ and which have an endpoint equal to $x$ (note that both endpoints of $\ell$ lie in
$\partial P$). Then $\delta_P(x)=\inf \{ | \frac{\ell}{\| \ell\| } \cdot H |\colon
\ell \in M_P(x) \wedge H\in \cH_P(x) \}$. Let $M_P=\bigcup_{x \in \partial P} M_P(x)$, and
$\cH_P=\bigcup_{x \in \partial P} \cH_P(x)$. Let $\delta_P=\inf \{ \delta_P(x)
\colon x \in \partial P\}$. 
\end{definition}

We topologize $M_P$ using the Hausdorff metric on the set of line segments in $\R^d$.
Similarly, we topologize $\cH_P(x)$ using the Hausdorff metric on the unit inward (toward $P$)  normal
vectors for the hyperplanes.

\begin{lemma} \label{lem:mlc}
$M_P$ is a closed set, and in particular $M_P(x)$ is closed for all $x \in \partial P$.
\end{lemma}

\begin{proof}
It is enough to show that if $\ell_n=\ell_{x_n,y_n}$ are line segments in $M_P$
and $\ell_n \to \ell=\ell_{x,y}$, then $\ell \in M_P$. From our remark above that the lengths
of the line segments in $M_P$ are bounded below, we have that the length of $\ell$
is positive, that is, $x \neq y$. Suppose $\ell \notin M_P$. Let $\ell'= \ell_{x',y'}=
z' + s' \ell \subseteq P$ where $s'>1$. Let $w \in \interior{P}$.
Let $B\subseteq P$ be an open ball containing $w$. 
By convexity, the convex hull $C$ of $B$ and $\ell'$ is contained in $P$.
But then $\interior{C}$ contains a translated (not enlarged) copy of $\ell$. For $n$ large enough $\ell_n$
will then also lie in $\interior{C}$, a contradiction to $\ell_n$ being maximal. 
\end{proof}

\begin{lemma}
For any compact, convex set $P\subseteq \R^d$, and any point $x \in \partial P$,
there is a line segment $\ell \in M_P(x)$ and a supporting hyperplane
$H\in \cH_P(x)$ such that $\delta_P(x)= | \frac{\ell}{\| \ell \| } \cdot H |$.
\end{lemma}

\begin{proof}
This follows immediately from the fact that $M_P(x)$ and $\cH_P(x)$ are compact
and the map $(\ell, H) \mapsto \frac{\ell}{\| \ell \| } \cdot H$ is continuous. 
\end{proof}

\begin{lemma} \label{lemma:delta}
Let $P$ be a compact, non-degenerate, strictly convex set. Then
$\delta_P >0$.
\end{lemma}

\begin{proof}
Suppose $\delta_P=0$. Let $\ell_n=\ell_{x_n,y_n}\subseteq P$ be maximal line segments with
respect to $P$, and $H_n \in \cH_P(x_n)$ be supporting hyperplanes with
$\frac{\ell_n}{\| \ell_n\| } \cdot H_n \to 0$. 
Without loss of generality we may assume that $\ell_n \to \ell=\ell_{x,y}$ (with $x \neq y$) and
the $H_n$ converge to $H$, a supporting hyperplane for $P$ through the point $x$. 
By Lemma~\ref{lem:mlc}, $\ell \in M_P(x)$. Also, $H\in \cH_P(x)$ and so by continuity of the
dot product we have that $\frac{\ell}{\| \ell \| }\cdot H=0$. So, $\ell \subseteq H\cap P$,
but $H\cap P\subseteq \partial P$, and so $\ell \subseteq \partial P$. This contradicts the
strict convexity of $P$. 
\end{proof}

Suppose now $P\subseteq \R^d$ is a compact, strictly convex set. Consider the no-$\beta$ McMullen
game for target set $A \subseteq \R^d$ with respect to the set $P$. If player $\I$ wins
this game, then by Fact~\ref{psf} we have that $A$ contains a perfect set.
For the rest of the argument we assume that $\II$ has a winning strategy $\tau$. 
We must show that $A$ is countable. 
Suppose towards a contradiction that $A$ is uncountable. 
Let
\[
m=\min\setof{k \in \omega \suchthat \exists H~\text{an affine subspace},~\dim{H}=k,~\text{and}~A \cap
  H~\text{is uncountable}}.
\]
Let $H$ be a $m$-dimensional affine subspace of $\mathbb{R}^d$ so that
$H \cap A$ is uncountable.  Without loss of generality, we can assume
that $A \subseteq H$.

Next consider the tree $T$ of all even length positions in the game
which are consistent with $\tau$, and for each $x \in A$, let $T_x$ be
the subtree of $T$ consisting only of positions $p=(B_0, B_1, \dots
B_{n+1})$ so that $x \in \interior {B}_n \setminus \closure{B_{n+1}}$,
i.e., the point $x$ is well-inside player $\I$'s last move, and was not
yet deleted by $\tau$.  Since $\tau$ is a winning strategy for $\II$,
no full infinite run of the game can exist which results in a point $x
\in A$, and thus for each $x \in A$, the tree $T_x$ is wellfounded.

Let $N$ be a countable elementary substructure containing $\tau$ and
$H$, and let $T^N=T \cap N$ and for each $x \in A$ let $T_x^N = T_x \cap N$.
Note that $T$ is definable from $N$ and so $T\in N$, but $T_x \notin N$ in general.
The set of all
positions in $N$ is countable, and so the set of all positions which
occur as terminal nodes in some tree $T_x^N$ for $x \in A$ is also
countable.  Since $A$ is uncountable by hypothesis, we can find some
fixed position $p=(B_0, \dots, B_{n+1})$ and some uncountable subset
$A' \subseteq A$ so that $p$ is terminal in $T_x^N$ for every $x \in
A'$.

Note that by the choice of $A'$, we have that $A' \subseteq
\interior{B}_n$, and since $A'$ is uncountable, it must have a strong
limit point $x_0 \in A'$.  Since $x_0 \in \interior{B}_n \setminus
\closure{B_{n+1}}$, there is some small enough $\epsilon$ so that
$B(x_0, \epsilon) \subseteq \interior{B}_n \setminus
\closure{B_{n+1}}$.  Next, we note that since $A' \subseteq A
\subseteq H$, and $x_0$ is a strong limit point of $A'$, we know that
$B(x_0, \epsilon) \cap H \cap A'$ is uncountable, and $p \concat
B(x_0, \epsilon)$ is a legal position in the game.  Let $D$ be a small
enough rational ball containing $x_0$ so that for any $x, y \in D$,
$[x, y]  \subseteq B(x_0, \epsilon)$. Note that $\tau$ induces a one-round
strategy in the no-$\beta$ McMullen game on $D$.

We will need the following two lemmas
to complete the proof. For the first lemma, we introduce the following terminology.

Given hyperplanes $H_1,\dots,H_k$ in $\R^d$, we define the family of {\em cones} which they induce
as follows. For each hyperplane $H$, we let $H^+$ and $H^-$ be the two closed half-spaces in $\R^d$
determined by $H$. For each $p \in \{ \pm 1\}^k$, let $C(p)= \bigcap_i H^{p(i)}$, where by $H^{+1}$
be mean $H^+$ and by $H^{-1}$ we mean $H^{-}$. We let $-p$ be the sequence $(-p)(i)=-p(i)$.
We let $C^{\pm}(p)$ be the union of $C(p)$ and $C(-p)$. We refer to $C^{\pm}(p)$
as a cone determined by the hyperplanes $H_1,\dots, H_k$. We arbitrarily choose from among
$p $ and $-p$ and let $C^+(p)=C(p)$ be the ``positive cone'' associated to the $H_i$, and
$C^-(p)=C(-p)$ be the ``negative cone'' associated to the $H_i$.  In this manner every collection
$H_1,\dots,H_k$ of hyperplanes induces a collection $C_1,\dots,C_\ell$ of cones in $\R^d$, and
for each of these cones $C_i$ we have defined the positive part $C_i^+$ and negative part $C^{-}_i$.
The following ``cone lemma'' will be proved in the next section.

\begin{lemma} \label{cone}
Suppose $A \subseteq \mathbb{R}^d$ is uncountable.  Suppose that $H_1,
\dots H_k$ are distinct hyperplanes inducing the family of cones $C_1
\dots C_\ell$.  Then there is a cone $C_i$ and an uncountable subset
$A' \subseteq A$ so that for every point $x \in A'$, $x$ is a strong limit
both from $C_i^+(x) \cap A'$ and $C_i^-(x) \cap A'$.
\end{lemma}

Using Lemma~\ref{cone} we next prove the following lemma concerning elementary
substructures and their relation to the no-$\beta$ McMullen game.

\begin{lemma}\label{lemma:cantdelunc}
Suppose $\tau$ is a one-round strategy in the no-$\beta$
McMullen game for player $\II$ on the ball $B \subseteq
\mathbb{R}^d$.  Then for any elementary substructure $N$ containing
$\tau$ and $P$, and for any uncountable set $A \subseteq B$ which contains at
most countably many points from any affine hyperplane, there is some
legal move (for $\I$) $Q=w+s P \subseteq B$ from $N$ so that $A \cap Q
\setminus \tau(Q)$ is nonempty.
\end{lemma}

\begin{proof}[Proof of Lemma~\ref{lemma:cantdelunc}]
Let $D\subseteq B$ be a rational ball so that $A \cap D$ is uncountable and
for every $x,y \in D$, $\bp{x}{y} \subseteq B$. 
Recall $P$ is compact and strictly convex. Let $\delta=\delta_P>0$ be as in Lemma~\ref{lemma:delta}.
Let $H_1,\dots,H_k$ be distinct hyperplanes in $\R^d$ with corresponding
cones $C_1,\dots,C_\ell$, such that for any of the cones $C_i$ and
any two unit vectors $u,v \in C_i$, $1-| u \cdot v|  < \frac{1}{2}\delta^2$. This is easily done by
taking sufficiently many of the hyperplanes. Fix the cone $C=C_i$ from Lemma~\ref{cone},
along with the uncountable set $A' \subseteq A$. 
Let $x \in A' \cap D$. Let $y \in A' \cap D$ be such that if $v$
is the vector from $x$ to $y$ then $v \in C^{\pm}$. 
Let $R= \bp{x}{y}$, so $R\subseteq B$.  Let $R'=\tau(R)$. Since $R'$ is a smaller scaled copy than $R$,
it follows from the definition of $\bp{x}{y}$ that at least one of $x,y$ is not in $R'$. The rest of the argument
is symmetrical between $x$ and $y$, so we assume that $x \in R \sm R'$.

Say without loss of generality
that $v \in C^+$. By definition of $\bp{x}{y}$, $\ell_{x,y} \in M_R(x)$. Since $\delta_P(x) \geq \delta_P =\delta$,
it follows from the definition of $\delta_P(x)$ that 
for any supporting hyperplane $H$ for $R$ at $x$,
we have that $|\frac{\ell}{\|\ell\|} \cdot H| > \delta$. 
This implies that for any point $w \in C^+(x)$ that if $\ell'=\ell_{x,w}$, then 
$|\frac{\ell'}{\|\ell'\|} \cdot H| > 0$. This is because $| \frac{\ell}{\|\ell\|} \cdot 
\frac{\ell'}{\|\ell'\| }|>1-\frac{1}{2}\delta^2$, $| \frac{\ell}{\|\ell\|} \cdot H| >\delta$,
and $| \frac{\ell'}{\|\ell'\|} \cdot H|\geq | \frac{\ell}{\|\ell\|} \cdot H| -
\| \frac{\ell}{\|\ell\|}- \frac{\ell'}{\|\ell'\|} \|$, and 
$\| \frac{\ell}{\|\ell\|}- \frac{\ell'}{\|\ell'\|} \|^2 =2(1-\frac{\ell}{\|\ell\|} \cdot 
\frac{\ell'}{\|\ell'\| })<\delta^2$. 
So, $|\frac{\ell'}{\|\ell'\|} \cdot H|>0$.

In particular, for any point
$z \in C^+(x)$ which is sufficiently close to $x$ we have that $z \in \interior{R}$. 
To see this, we first observe that for every $z \in C^+(x)$ there is a $z' \in \ell_{x,z}$
such that $\ell_{x,z'}\subseteq R$. If not, then $\ell_{x,z}\cap R=\{ x\}$. There is then 
a supporting hyperplane $H$ for $R$ which contain $x$ and $\ell_{x,z} \sm \{x\}$ 
lies (non-strictly) on the other side of $H$ from $R$.  Since $y,z \in C^+(x)$ and 
$C^+(x)$ is convex, the line segment $\ell_{y,z}$ lies in $C^+(x)$. Thus
there is a point $w \in \ell_{y,z}\subseteq C^+(x)$ with $w \in H$. This contradicts 
the above fact. Next we observe that there is a $z'' \in \ell_{x,z}$
with $\ell_{x,z''}\subseteq \interior{R}$. It is enough to see that there is a 
$z'' \in \ell_{x,z'}$ with $z'' \in \interior{R}$. If this failed, then 
$\ell_{x,z'}$ is contained in $\partial R$.  Then there is a supporting hyperplane 
$H$ for $R$ which contains $\ell_{x,z'}$, a contradiction as $z' \in C^+(x)$.

Fix a rational ball $B_1$ about $x$ so that $R' \cap B_1=\emptyset$. Let $z \in B_1\cap C^+(x)\cap A'$,
with $z \in \interior{R}$, which we can do
by the choice of $C$ as $x \in A'$. Let $B_2 \subseteq B_1$ be a small enough rational ball about $z$
so that $B_2 \subseteq \bp{x}{y}=R$. See Figure~\ref{lemfig} for an illustration of these objects.

\begin{figure}
  \begin{tikzpicture}[scale=4]
    \draw (0,0) circle (0.75 cm);
    \draw [fill=gray!20!white] (-0.05,0) ellipse (0.75*0.48cm and 0.75*0.24cm);
    \draw [fill=black, draw=black] (-0.5, -0.3) circle (0.015cm);
    \node [anchor=north] at (-0.5, -0.3) {$x$};
    \draw [dashed, draw=black!60!white] (-0.5 -0.25, -0.3 -0.25) -- (-0.5 +1, -0.3 +1); 
    \draw [dashed, draw=black!60!white] (-0.5+1*-0.4, -0.3+0.35*-0.4) -- (-0.5 +1*1.5, -0.3 +0.35*1.5);
    \draw [fill=black, draw=black] (-0.5+0.65, -0.3+0.35) circle (0.015cm);
    \node [anchor=north] at (-0.5+0.65, -0.3+0.35) {$y$};
    \draw (-0.18,-0.125) ellipse (0.48cm and 0.24cm);
    \node [anchor=north] at (0, -0.35) {$\bp{x}{y}$};
    \draw [fill=black, draw=black] (-0.5+0.15, -0.3+0.1) circle (0.015cm);
    \draw [thin] (-0.43, -0.29) circle (0.18cm);
    \draw [thin] (-0.36, -0.22) circle (0.06cm);
  \end{tikzpicture}
  \caption{Proof of Lemma~\ref{lemma:cantdelunc}} \label{lemfig}
\end{figure}
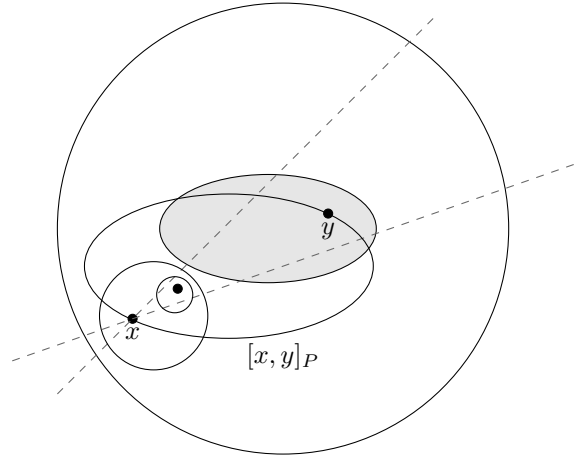

Consider the statement
\[
\varphi(x,y)=
(x,y \in D) \wedge (\tau(\bp{x}{y}) \cap B_1=\emptyset) \wedge (B_2 \subseteq \bp{x}{y}).
\]
By elementarity of $N$
we get that
\[
\exists x', y' \in N\ (x',y' \in D)\wedge (\tau(\bp{x'}{y'}) \cap B_1=\emptyset) \wedge (B_2 \subseteq \bp{x'}{y'}).
\]
Let $Q=\bp{x'}{y'}$. Then we have $z \in \bp{x'}{y'}$ since $z \in B_2 \subseteq \bp{x'}{y'}$ and also
$z \notin \tau(\bp{x'}{y'})$ since $z \in B_2 \subseteq B_1$ and $\tau(\bp{x'}{y'})\cap B_1=\emptyset$. 
Also, $x',y'\in D$ so $Q=\bp{x'}{y'}\subseteq B$. Since $x',y' \in N$, $Q \in N$. As $z \in A' \cap (Q\sm \tau(Q))$,
this proves the lemma. 

\end{proof}

By Lemma~\ref{lemma:cantdelunc}, we know that there is some move $Q
\subseteq D$, $Q\in N$,  so that $A' \cap Q \setminus \tau (Q)$ is nonempty.
This contradicts that the
position $p$ was terminal for every $x \in A'$.
This completes the proof of Theorem~\ref{thma}.
\\

\hfill $\square$

\section{Cone Lemma}

In this section we first prove the cone lemma, Lemma~\ref{cone}, which for convenience
we restate. We then prove another variation of the theorem which abstracts part of the argument
and does not mention cones. We will make use of this other variation as well in the proof 
of Theorem~\ref{thmb}.
\begin{lemma*} 
Suppose $A \subseteq \mathbb{R}^d$ is uncountable.  Suppose that $H_1,
\dots H_k$ are distinct hyperplanes inducing the family of cones $C_1
\dots C_\ell$.  Then there is a cone $C_i$ and an uncountable subset
$A' \subseteq A$ so that for every point $x \in A'$, $x$ is a strong limit
both from $C_i^+(x) \cap A'$ and $C_i^-(x) \cap A'$.
\end{lemma*}

\begin{proof}
Suppose $A \subseteq \R^d$ is uncountable, and the hyperplanes $H_1,\dots,H_k$
and given which induce the cones $C_1, \dots, C_\ell$ in $\R^d$.

Without loss of generality we may assume that $A\cap H$ is countable for
every affine hyperplane $H \subseteq \R^d$. To see this, let $S\subseteq \R^d$
be an affine subspace of $\R^d$ of minimal dimension such that $A\cap S$ is
uncountable. By an affine translation we may identify $S$ with $\R^m$
for some $m =\dim(S)<d$. This maps the hyperplanes $H_1,\dots,H_k$
to new hyperplanes $H'_1,\dots,H'_k$ and gives new cones $C'_1,\dots,C'_\ell$.
By the same affine translation the set $A$ moves to the uncountable set $A'$.
$A'$ now intersects every hyperplane in $\R^m$ in a countable set.
Proving the result for $A'$ then easily gives the result for $A$ by moving back
under the affine translation. So, we assume for the rest of the argument that
$A\cap H$ is countable for every affine hyperplane $H\subseteq \R^d$.

Let $A_0=A$.
Consider the first cone $C_1$. We ask if there is an uncountable $A'_0 \subseteq A_0$
such that for every $x \in A'_0$ there is a neighborhood $U$ of $x$ such that
either $U\cap C^+_1(x)\cap A'_0$ is countable or $U\cap C^-_1(x)\cap A'_0$ is countable.
If this is the case then we may
thin $A'_0$ out to an uncountable set $A_1 \subseteq A'_0$ and fix an integer
$p(1)\in \{ \pm 1\}$ so that for all $x \in A_1$, there is a neighborhood $U$ of $x$
such that $U \cap C^{p(1)}_1(x) \cap A_1$ is countable. That is, by thinning from $A'_0$ to $A_1$ we fix the
side of the cone $C_1$ which is weakly isolating $x$ in $A_1$ with respect to $C_1$.
We continue in this manner for $\ell$ stages, assuming that we can at each stage
go from $A_k$ to $A'_k\subseteq A_k$ using the cone $C_k$. We then define $A_{k+1} \subseteq A'_k$
and $p(k)$ exactly as in the first step. Assuming this construction continues
for all $\ell$ steps, we end with an uncountable set $A_\ell \subseteq A$ and a
``cone pattern'' $p \in \{ \pm 1\}^\ell$.

In this case, for each $x \in A_\ell$ we let $B_x$ be a rational ball containing $x$
such that $\forall i\ B_x \cap C^{p(i)}_i(x) \cap A_\ell$ is countable. We may fix a rational
ball $B$ so that $B_x=B$ for an uncountable subset $A'_\ell$ of $A_\ell$.
Let $x,y \in A'_\ell$ be distinct strong limit points of $A'_\ell$, so
in particular $x,y \in B=B_x=B_y$.
Since the set of points of $A'_\ell$ which are strong limit points of
$A'_\ell$ is uncountable, from our assumption on $A$ we may assume that
the line $\ell_{x,y}$ between $x$ and $y$ does not lie on any of the $H_1,\dots,H_k$.
For one of the cones,
say $C_i$, we have that $y \in C_i^{\pm}(x)$ (and so also $x \in C^{\pm}_i(y)$). Say without loss of
generality $y \in C^+_i(x)$, and so $x \in C^-_i(y)$. If $p(i)=+1$, then we have a contradiction
from $y \in C^+_i(x)$, and so $y$ is in the interior of $C^+_i(x)$, 
and the facts that $x \in A'_\ell$ and $y$ is a strong limit point of $A'_\ell$.
Similarly, if $p(i)=-1$ we
have a contradiction from the fact that $x \in C^-_i(y)$ and $y \in A'_\ell$.

Thus, the above construction cannot proceed through all $\ell$ cones. Let $i$
be the least stage where the construction fails. This means that there is an
uncountable set $A'=A_{i-1}\subseteq A$ such that for any uncountable $A'' \subseteq A'$
there is an $x \in A''$ such that $x$ is a strong limit of both $C^+_i(x) \cap A''$ and
$C^-_i(x)\cap A''$. Let $E\subseteq A'$ be the set of $x \in A'$ such that there is a neighborhood $U$
of $x$ such that either $U \cap C^+_i(x)\cap A'$ is countable or $U \cap C^-_i(x)\cap A'$ is countable.
$E$ is countable, as otherwise we could fix a rational ball $B$ and an uncountable $E' \subseteq E\cap B$
such that for all $x \in E'$, either $C^+_i(x)\cap E'$ is countable or $C^-_i(x)\cap E'$ is countable.
This contradicts that the construction failed at step $i$,
as we could take $A_i=E'$.
So, $A'\sm E$ is uncountable. From the definition of
$A'$ and the fact that $E$ is countable we also have that for any $x \in A'\sm E$, $x$
is a strong limit point of both $(A'\sm E) \cap C^+_i(x)$ and $(A'\sm E) \cap C^-_i(x)$ (note that being a
strong limit point of $A'\cap C^+_i$ implies being a strong limit point of $(A'\sm E)\cap C^+_i$
as $E$ is countable). So, the set $A'\sm E$ and the cone $C_i$ verify the statement of the theorem.

\end{proof}

The next theorem is the abstraction of Lemma~\ref{cone}. In the statement
of the theorem, by $A^2$ we mean the set of unordered pairs from the set $A$.

\begin{theorem} \label{gencone}
Let $X$ be a second countable space, $A\subseteq X$ uncountable, and $c \colon A^2 \to \{ 1,\dots, n\}$
a partition of the (unordered) pairs from $A$ into finitely many colors. Then there is an uncountable
$B\subseteq A$ and an $i \leq n$ such that $B$ is partially homogeneous for $i$. By this we mean that 
for every $x \in B$, $x$ is a strong limit point of the set $\{ y\in B\colon c(x,y)=i\}$. 
\end{theorem}

\begin{proof}
The proof is similar to that of Lemma~\ref{cone}, so we only give a sketch. 
Starting with $A_0=A$, and given at stage $i$ the sets $A_{i-1} \subseteq A_{i-2}\subseteq \cdots \subseteq A_0$, 
we ask if there is an uncountable $A_i \subseteq A_{i-1}$ such that 
for every $x \in A_i$ there is a neighborhood $U$ of $x$ such that 
$\{ y \in U \cap A_i\colon c(x,y)=i\}$ is countable. If we do this for all $i \leq n$,
let $A_n$ be the ending set, which is uncountable and such that for every $i \leq n$ 
and every $x \in A_n$, there is a neighborhood $U$ of $x$ 
such that the set $\{ y \in A_n \colon c(x,y)=i\}$  is countable. Since $X$ is second countable
we may fix the neighborhood $U$ which has this property for every point in an uncountable subset $A'_n \subseteq A_n$. 
Fix an $x \in A'_n$. Then for each $i \leq n$ there are only countably many $y \in U\cap A'_n$ with $c(x,y)=i$,
and thus only countably many points of $U\cap A'_n$, a contradiction.

So, the construction must stop with some set $A_{\ell-1}$ (that is, $A_\ell$ is not defined), for some $\ell\leq n$. 
Let $C$ be the set of $x \in A_{\ell-1}$ such that there is a neighborhood $U$ of $x$ with 
$\{ y \in A_{\ell-1}\cap U\colon c(x,y)=\ell\}$ is countable. We have that $C$ is countable,
otherwise we could use $C=A_\ell$ and continue the construction. 
So, $C$ is countable and if we let $B=A_{\ell-1} \sm C$, then easily $B$ is a partially homogeneous set for $\ell$. 
\end{proof}

\section{Polytopes in $\R^d$} \label{sec:poly}

We begin this section with some notation and statement of some results concerning 
polytopes which we need for the theorem. These geometric results on polytopes will be proved in the next section.

\begin{definition}
Let $x, y \in \partial P$ be such that $\ell_{x, y}$ is maximal.  
We say $\ell_{x,y}$ is an \emph{extreme}  line segment if 
for every $t \in \R^d$ and every $\epsilon >0$, either $\ell_{x,y} +\epsilon t$ 
or $\ell_{x,y} -\epsilon t$ does not lie in $P$.
\end{definition}

\begin{definition} \label{Gdef}
Let $\mathfrak{F}$ be a subset of the faces of $P$.
We let $G(x,\mathfrak{F})$ be the set 
of points $y$ such that there is a scaled, translated copy $P'$
of $P$ with $\ell_{x,y}$ extreme and maximal in $P'$ and for every face $F'$ of $P'$, $\ell_{x,y}\subseteq F'$ iff 
the corresponding face $F$ of $P$ is in $\mathfrak{F}$. 
\end{definition}

Given an affine subspace $S\subseteq \R^d$ and an intersection $\mathcal{F}$ of faces of $P$, 
we say $S$ is {\em parallel} to $\mathcal{F}$ if a translation of $S$ is contained in the 
affine extension of $\mathcal{F}$. By the {\em affine extension} of a subset 
$\mathcal{F}$ we mean the set of points of the form $x+t(y-x)$
where $x,y \in \mathcal{F}$ and $t \in \R$.

The next lemma is a generalization
of the cone lemma, Lemma~\ref{cone},  which we need for the main proof.
The proof will be  given in the next section.

\begin{lemma} \label{conepart}
Let $P$ be a polytope in $\R^d$. Suppose $A \subseteq S$ is uncountable where $S\subseteq \R^d$ is an affine subspace,
and assume that every hyperplane in $S$ intersects $A$ in a countable set.
Suppose that $H_1, \dots H_k$ are distinct hyperplanes in $S$ inducing the family of cones $C_1
\dots C_\ell$. 
Then there is a cone $C_i$, a subset $\mathfrak{F}$ of the faces of $P$  and an uncountable subset
$A' \subseteq A$ so that for every uncountable $B\subseteq A'$ there is an uncountable $D\subseteq B$
such that for every point $x \in D$, $x$ is a strong limit
both from $C_i^+(x) \cap D$ and $C_i^-(x) \cap D$ of points $y$ in $G(x,\mathfrak{F})$.
\end{lemma}

In \S\ref{sec:sc} we gave the proof of the main theorem in the case where $P$ is a compact, strictly convex
set in $\R^d$, Theorem~\ref{thma}. In this section we modify the argument to handle the case where $P$
is a $d$-polytope in $\R^d$. By a $d$-polytope (generalized polyhedron) we mean a compact set $P$ in $\R^d$ which is the
intersection of a finite set of closed half-spaces, and which has non-empty interior.
Such a set is necessarily convex, but not strictly convex.
Equivalently, $P$ is the convex hull of a finite set in $\R^d$
which does not lie in a $<d$ dimensional affine subspace. 
For the rest of this section $P$
will denote a fixed $d$-polytope in $\R^d$, and we prove our main theorem, Theorem~\ref{thmb}, in this case.

As in the proof of Theorem~\ref{thma} we assume $A\subseteq \R^d$ is uncountable and $\tau$ is a winning
strategy for $\II$ in the no-$\beta$ McMullen game and we proceed to get a contradiction.
Fix an affine subspace  $S \subseteq \R^d$ of minimal dimension 
such that $A\cap S$ is uncountable. So for every affine hyperplane $H \subseteq \R^d$, either 
$S\subseteq H$ or $H\cap A$ is countable. 
Replacing $A$ by $A\cap S$, we may assume that $A\subseteq S$.

Let $H_1,\dots, H_k$ be the affine hyperplanes determined by the faces of $P$ such that a translation of $H_i$ 
intersect $S$ properly. Applying Lemma~\ref{conepart} gives a set of faces $\mathfrak{F}$, 
a cone $C$ defined by the $H_i$, and an uncountable $A'\subseteq A$. To save notation,
let us rename the set $A'$ to $A$.

For $x,y \in A$ let $[x,y]_P$ be a scaled, translated copy of $P$ for which $\ell_{x,y}$
is maximal, and if $y\in G(x,\mathfrak{F})$ (equivalently, $x \in G(y,\mathfrak{F})$), then we take $[x,y]_P$ so that 
$\ell_{x,y}$ is on exactly the faces in $\mathfrak{F}$.

We define the trees $T$ and $T_x$ as in Theorem~\ref{thma}, with a few changes. 
As before, $p=(P_0,\dots, P_{n+1})$ of even length is in $T$ iff it is consistent with $\tau$, and
we also require that every move $P_{2i}$ is of the form $\bp{x}{y}$ for some $x,y \in A$ with $y \in G(x,\mathfrak{F})$.
Such a $p \in T$ is in $T_x$ if in addition 
$x$ is in the relative interior $\interior{P}_n(S)$ 
of $P_n \cap S$ with respect to the relative 
topology on $S$, and $x$ is not in $P_{n+1}$. 
So again we have that for every $x \in A$
the tree $T_x$ is wellfounded.

We let $N$ be a countable elementary substructure containing $\tau$,
$P$, $A$, $S$, and the map $x,y\mapsto \bp{x}{y}$.  
As before let $T^N=T\cap N$ and $T^N_x=T_x \cap N$. Also as before we can find a position
$p=(P_0,P_1,\dots, P_{n+1})$ (of even length; recall all of the
moves $P_i$ are translated, scaled copies of $P$) such that for some uncountable $A' \subseteq A$
and all $x \in A'$, $p$ is a terminal position in $T^N_x$. In particular, for all $x \in A'$ we have
$x \in \interior{P}_n(S) \sm P_{n+1}$.  From Lemma~\ref{conepart} we may assume that $A'$
has the property that every $x \in A'$ is a strong limit point from both of
$A' \cap C^+(x)\cap G(x,\mathfrak{F})$ and $A'\cap C^-(x)\cap G(x,\mathfrak{F})$.

Let $D$ be a rational ball in $\R^d$ such that 
$D\cap S \subseteq {P_n^\circ(S) }\sm P_{n+1}$, 
$A'\cap D$ is uncountable. We may take $D$ to be sufficiently small so that 
for any distinct $x,y \in D\cap S$, if $y\in G(x,\mathfrak{F})$ 
(so $\ell_{x,y}$ lies on exactly the faces in $\mathfrak{F}$), 
then $\bp{x}{y}$ is disjoint from any face of $P_{n+1}$ and any face of $P_n$
not in $\mathfrak{F}$. To see this, let $w$ be a strong limit point of $A'$ in 
$\interior{P}_n(S) \sm P_{n+1}$. Let $D$ be a small enough neighborhood of $w$ 
such that $D$ is disjoint from $P_{n+1}$ and also disjoint from any face 
$F$ of $P_n$ which $w$ does not lie on. If there is a face $F$ of $P_n$ 
not in $\mathfrak{F}$ which $w$ lies on, then the affine extension $H$ of $F$ must contain 
$S$ as otherwise its intersection with $S$ is proper, and since $w \in \interior{P}_n(S)$ 
we would have $w \notin F$. If $P_n=\bp{\bar x}{\bar y}$ with $\bar y \in G(\bar{x}, \mathfrak{F})$, 
then $\ell_{\bar{x}, \bar{y}} \nsubseteq F$ and so also $\ell_{\bar{x}, \bar{y}}$ is not a subset of 
$H$. Since $\bar{x}$, $\bar{y}\in S$, this shows that $H$ properly intersects $S$, a contradiction.

It follows that for any $x, y \in A'\cap D$, 
if $y \in G(x,\mathfrak{F})$ then $\bp{x}{y}$ is disjoint from
all faces of $P_{n+1}$ and those of $P_n$ not in $\mathfrak{F}$, 
and for those faces $F\in \mathfrak{F}$ we have that the affine extension of the face $F'$ of 
$\bp{x}{y}$ corresponding to $F$ is the affine extension of the face $F$ of $P$
(note that the affine extension $H$ of $F$ and the affine extension $H'$ of $F'$ cannot
be disjoint as $H$ contains $S$ and $x,y \in S\cap F'$). 
So, for any such $x,y \in A'\cap D$ (with $y \in G(x,\mathfrak{F})$) 
we have that $\bp{x}{y} \subseteq P_n\sm P_{n+1}$. Thus, $\bp{x}{y}$
is a valid next move for $\I$ in the game $G_P(A)$.

Fix $x,y \in A'\cap D$ with $y \in G(x,\mathfrak{F})$. This is possible since 
$A'\cap D$ is uncountable and each of the co-dimension one affine hyperplanes $H(x)$
from Lemma~\ref{dimleme} has a countable intersection with $A'$. Again, 
$R=\bp{x}{y}\subseteq P_n\sm P_{n+1}$ is a valid move for $I$. Let $R'=\tau(R)$ 
be $\II$'s response by $\tau$. By maximality of $\ell_{x,y}$ in $R$, 
$R'$ cannot delete both $x$ and $y$, so without loss of generality assume $x \in R'$. 
Also without loss of generality assume $y  \in C^+(x)$. 
Since $x \in A'$, $x$ is a strong limit point of points $z \in A'\cap  C^+(x)$. 
For such $z$ close enough to $x$ we have that $z \in R$. This is because if $F$ is a face of $R$
which does not contain $x$, and $z$ is close enough to $x$ then $z$ is on the same side of $F$ that $x$ is on. 
If $F$ is a face of $R$ containing $x$, and the affine extension of $F$ contains $S$
then since $z \in A'\subseteq S$, we have that $z$ is also on the same side of $F$ that $x$
is on. Finally, if $x \in F$ and the affine extension $H$ of $F$ does not contain $S$, then 
$H\cap S$ is one of the co-dimension one subsets of $S$ used in forming the cone $C^+(x)$, and 
so $z$ lies on the same side of $F$ as $R$ does (we use here that $y \in C^+(x)$). So, $z \in R$,
and we may assume $z \notin R'$ by taking $z$ close enough to $x$.

Let $B_1$ be a rational ball about $x$ so that $R' \cap B_1=\emptyset$.
We may assume that $z \in B_1$. Let $B_2 \subseteq B_1$
be a rational ball about $z$ with $B_2 \cap S \subseteq R$. 
Let $\varphi(x,y)$ be the statement
\begin{equation*}
\begin{split}
\varphi(x,y)& =
(x,y \in D \cap S) \wedge (y \in G(x,\mathfrak{F})) \wedge (B_2\cap S \subseteq \bp{x}{y})
\\ & \quad \wedge (\tau(\bp{x}{y}) \cap B_1=\emptyset) \wedge ( \bp{x}{y} \subseteq P_n\sm P_{n+1}).
\end{split}
\end{equation*}
Since $S, D, B_1, P_n, P_{n+1}=\tau(P_n)$, and the map $x,y\mapsto \bp{x}{y}$ are all in $N$, by elementarity we have that
\begin{equation*}
\begin{split}
\exists x',y' & \in N\cap D\cap S\ [ (y' \in G(x',\mathfrak{F})) \wedge (B_2\cap S \subseteq \bp{x'}{y'})
\\ & \quad \wedge (\tau(\bp{x'}{y'}) \cap B_1=\emptyset) \wedge ( \bp{x'}{y'} \subseteq P_n\sm P_{n+1})].
\end{split}
\end{equation*}
Let $P_{n+2}=\bp{x'}{y'}$. Since $\bp{x'}{y'} \subseteq P_n\sm P_{n+1}$,
$P_{n+2}$ is a legal move for $\I$ in $G_P(A)$ extending the position $p$.
Since $x',y' \in N$, we have that $P_{n+2}\in N$. Since $z \in B_2 \cap S$ and
$B_2 \cap S \subseteq \bp{x'}{y'}$, we have that $z \in P_{n+2}$. Since
$z \in B_1$ and $\tau(\bp{x'}{y'}) \cap B_1=\emptyset$ we have that $z \notin  \tau(\bp{x'}{y'})=\tau(P_{n+2})$.
As $P_{n+2}\in N$ and $z \in A'$, and since $y'\in G(x',\mathfrak{F})$, this contradicts the fact that $p$ was a terminal node
of $T^N_z$.

This completes the proof of Theorem~\ref{thmb}.

\section{Lemmas on Polytopes and a Generalized Cone Lemma} \label{sec:polylemmas}

\begin{notation}
Let $P$ be a $d$-polytope in $\mathbb{R}^d$.  For $x \in \partial P$, let $\mathcal{F}_x$ denote the 
intersection of all faces of $P$ which $x$ lies on.
\end{notation}

\begin{definition}
Let $x, y \in \partial P$ be such that $\ell_{x, y}$ is maximal.  
Let $\mathcal{F}$ be the intersection of a  (finite) subset of the faces of $P$.
We say $\ell_{x, y}$ is \emph{critical} for $\mathcal{F}$ of $P$ if there is a 
maximal line segment parallel to $\ell_{x, y}$ which lies in $\mathcal{F}$, and
for every neighborhoods $U$ of $x$ and $V$ of $y$ there are points
$u \in U$ and $v \in V$ such that some translated scaled copy of $\ell_{u,v}$ 
lies in $\mathcal{F}$ but no such copy is maximal in $P$.
\end{definition}

Note that a line segment being critical for  faces $\mathfrak{F}$ is invariant under
translation, provided the translated segment still lies in $P$.

\begin{lemma}\label{dimlem} 
Let $P \subseteq \mathbb{R}^d$ be a $d$-polytope, and
let $\ell_{x,y}$ be an extreme maximal line segment in $P$. 
Let $S \subseteq \mathbb{R}^d$ be an affine linear subspace.  
Then $\dim{(\mathcal{F}_x \cap S)} + \dim{(\mathcal{F}_y \cap S)} \leq \dim{S} -1$.
\end{lemma}

\begin{proof}
Suppose not, so that $\dim{(\mathcal{F}_x \cap S)} +
\dim{(\mathcal{F}_y \cap S)} \geq \dim{S}$.  Let $B_x, B_y$ be bases
for the linear subspaces corresponding to $\mathcal{F}_x \cap S$,
$\mathcal{F}_y \cap S$ respectively, and let $B$ be a basis for the
translate of $S$ which contains the origin.

First suppose $\lspan{B_x} \cap \lspan{B_y} \neq \{ \vec 0\}$.  If we
let $v \in \lspan{B_x} \cap \lspan{B_y}$, then for all small enough
$\epsilon>0$, $x\pm \epsilon v \in \mathcal{F}_x \cap S$ and $y\pm\epsilon
v \in \mathcal{F}_y \cap S$. In particular, $x\pm \epsilon v \subseteq P$ and 
$y \pm \epsilon t \subseteq P$. This contradicts $\ell_{x,y}$ being extreme.

So we must have $\lspan{B_x} \cap \lspan{B_y} = \{\vec 0\}$, but then
we have that $B_x \cup B_y$ is linearly independent and contains at
least as many vectors as $B$ does, and so we have $\lspan{(B_x \cup B_y)} =
\lspan{B}$.  In particular, since $x, y \in S$, we have $y-x \in
\lspan B = \lspan{(B_x \cup B_y)}$, and so there are vectors $u \in
\lspan{B_x}$, $v \in \lspan{B_y}$ so that $u+v = y-x$.  For any
$\epsilon>0$ small enough, we have that $x-\epsilon u \in
\mathcal{F}_x \cap S$ and $y+\epsilon v \in \mathcal{F}_y \cap S$, and
by convexity, we have that the line segment $\ell_{x-\epsilon u,
  y+\epsilon v}$ is in $P$, but note that $(y+\epsilon v)-(x-\epsilon
u) = y-x+\epsilon(u+v) = (1+\epsilon)(y-x)$, and so the line segment
$\ell_{x-\epsilon u, y+\epsilon v}$ is strictly longer than $\ell_{x,
  y}$, contradicting the maximality of $\ell_{x, y}$ in $P$.
\end{proof}

\begin{lemma} \label{translationa}
Let $\ell_{x,y}$ be maximal in $P$. Suppose $u_n\to x$ and $v_n \to y$. 
Then there are maximal line segments $\ell_{u'_n,v'_n}$ parallel to $\ell_{u_n,v_n}$
and a subsequence $i_n$ with $u'_{i_n}\to x$ and $v'_{i_n} \to y$. 
\end{lemma}

\begin{proof}
Let $\ell_{u'_n,v'_n}$ be maximal in $P$ and parallel to $\ell_{u_n,v_n}$. Let 
$u'_{i_n}, v'_{i_n}$ be a subsequence which converges to some $\ell_{x',y'}$ which is 
necessarily maximal in $P$ as well. Note that $\ell_{x',y'}$ is parallel to $\ell_{x,y}$, and in fact 
is of the same length as well by maximality of $\ell_{x,y}$ and $\ell_{x',y'}$. 
Let $\ell_{x,y}=\ell_{x',y'}+t$. By compactness, for any $\epsilon >0$ there is a $\delta >0$ such that 
if $u'_{i_n}, v'_{i_n}$ are in $B_\delta(x')$, $B'_\delta(y')$ respectively, then 
for any $\alpha \leq 1-\epsilon$ we have that if $u'_{i_n}+\alpha t$ is on a face $F$
of $P$, then so is $x'+\alpha t$, and likewise for $v'_{i_n}$ and $y'$. Since 
the $x'+\alpha t$ can be translated further in the direction of $t$ and remain in $P$, it follows 
that the same is true for $u'_{i_n}+\alpha t$, and like wise for $y'+\alpha t$ and $v'_{i_n}+\alpha t$. 
so, for all $\alpha \leq 1-\epsilon$, we have that $u'_{i_n}+\alpha t$, $v'_{i_n}+\alpha t$ are in $P$. 
In particular $u''=u'_{i_n}+(1-\epsilon) t, v''=v'_{i_n}+(1-\epsilon)t \in P$. Note that 
$\ell_{u'',v''}$ is also maximal in $P$. Letting $\epsilon$ go to $0$ now gives the result. 
\end{proof}

\begin{lemma} \label{translationb}
Suppose $\ell_{x,y}$ is an extreme maximal line segment. Let $U'$, $V'$ be neighborhoods of $x$ and $y$ respectively. 
Then there are neighborhoods $U \subseteq U'$, $V\subseteq V'$ of $x,y$ such that if $u\in U$, $v\in V$
and $\ell_{u,v}$ is maximal in $P$, then there is a translation $\ell_{u',v'}$ of $\ell_{u,v}$ with 
$u'\in U'$, $v'\in V'$ and $\ell_{u',v'}$ is extreme and maximal in $P$. 

\end{lemma}

\begin{proof}
Let $\epsilon>0$ be such that the $\epsilon$ neighborhoods of $x$ and $y$ are contained in 
$U'$, $V'$ respectively. 
Consider a translation vector $t \in \R^d$. For at least one of the points $x$, $y$ 
we have that this point cannot be translated in the $+t$ or the $-t$ direction and stay inside $P$. 
Without loss of generality, say $x$ cannot be translated in the $+t$ direction and stay inside of $P$. 
So, there is a face $F$ of $P$ such that $x \in F$ and for every $\eta>0$, $x+\eta t$ is on the side of the affine hyperplane $H$
defined by $F$ not containing $P$. Note that if $\vec {n}$ is the outward normal vector for $F$, then 
$t \cdot \vec{n} >0$. Thus there is a neighborhood $W_t$ of $t$ and a $\delta_t>0$ such that if 
$t'\in W_t$ then $x+ \frac{\epsilon}{2} t'$ is more than $\delta_t$ away from $H$ and on the 
side of $H$ not containing $P$. It follows that if $u$ is within distance $\delta_t$ of $x$ 
and $t' \in W_t$, then $u+ \frac{\epsilon}{2} t'$ is on the side of $H$ not containing $P$. 
This defines a cover $\{ W_t\}$ of the set of translation vectors with corresponding $\delta_t>0$.
By compactness, there is a $\delta>0$ such that for all $t$ we have that 
if $u \in B(x,\delta)$ and $v \in B(y,\delta)$  then at least one of $u\pm\frac{\epsilon}{2} t$, 
$v\pm\frac{\epsilon}{2} t$ is outside of $P$.

Suppose now $\ell_{u,v}$ is maximal in $P$ with $u \in B(x,\delta)$, $v \in B(y,\delta)$. 
If $\ell_{u,v}$ is not extreme, there is a vector $t$ such that $\ell_{u,v}\pm\eta t $ is in $P$ 
for some $\eta>0$. By the definition of $\delta$ we have, without loss of generality, that $\ell_{u,v}+\frac{\epsilon}{2} t$
is not contained in $P$. So there is a $\eta< \eta' <\frac{\epsilon}{2}$ such that $\ell_{u,v}+\eta' t$ 
is contained in $P$ and the total number of faces the endpoints lie on has strictly increased 
(note that since $\ell_{u,v}$ is maximal, if  $\ell_{u,v}\pm\eta t \subseteq P$ then $u\pm \eta t$, $v \pm \eta t$
lie on at least the same faces of $P$ that $u$, $v$ respectively lie on).

Let $\delta_0=\frac{\epsilon}{2}$. We have shown that there is a $\delta_1>0$ such that if 
$u \in B(x,\delta_1)$, $v \in B(y,\delta_1)$ and $\ell_{u,v}$ is maximal, then either $\ell_{u,v}$ is extreme there is a translation
$\ell_{u_1,v_1}$ of $\ell_{u,v}$ in $P$ such that the total number of faces the endpoints lie on has strictly increased 
Let $N$ be the number of faces of $P$. Repeating this argument at most $N$ times gives $\frac{\epsilon}{2}=\delta_0>\delta_1>
\cdots \delta_N>0$ such that if $\ell_{u,v}$ is maximal in $P$ with $u \in B(x,\delta_N)$, $v \in B(y,\delta_N)$,
then there is a translation $\ell_{u',v'}$ of $\ell_{x,y}$ in $P$ which is extreme and with 
$u'\in B(x,\frac{\epsilon}{2}$), $v'\in B(y,\frac{\epsilon}{2})$.

\end{proof}

\begin{lemma}\label{dimlemc}
Let $\mathcal{F}$ be an intersection of faces of the $d$-polytope $P \subseteq \mathbb{R}^d$.  
Let $S$ be a proper affine linear subspace of $\mathbb{R}^d$ which is parallel to $\mathcal{F}$.
Suppose $x, y \in S \cap \mathcal{F}$ are such that $\ell_{x, y}$ is extreme and critical for $\mathcal{F}$. 
Then $\dim{(\mathcal{F}_x \cap S)}+\dim{(\mathcal{F}_y \cap S)} \leq \dim{S}-2$.
\end{lemma}

\begin{proof}
Let $u_n \to x$,  $v_n \to y$ witness the criticality of $\ell_{x,y}$. 
In particular, there is a scaled, translated copy of $\ell_{u_n,v_n}$ which lies in $\mathcal{F}$,
but no such maximal copy of $\ell_{u_n,v_n}$ lies in $\mathcal{F}$.  
Let $u'_n, v'_n \in
\partial P$ be such that $\ell_{u'_n, v'_n}$ is parallel to
$\ell_{u_n, v_n}$ and is maximal in $P$. From Lemma~\ref{translationa} we may assume that 
$u'_x\to x$ and $v'_n \to y$. From Lemma~\ref{translationb} we may assume that the $\ell_{u'_n,v'_n}$
are extreme and maximal. Renaming the points, we may assume that $u_n\to x$, $v_n\to y$, 
$\ell_{u_n,v_n}$ is extreme and maximal, and by criticality that $\ell_{u_n,v_n}$
does not lie in $\mathcal{F}$. Since 
there is a scaled, translated copy of $\ell_{u_n,v_n}$ which lies in $\mathcal{F}$,
neither endpoint $u_n, v_n$ lies in $\mathcal{F}$.

Let $B$ be a basis for the translate of $S$ which contains the origin.
Define $S'=x+\lspan{(B \cup \{u_0 - x\})} = y+\lspan{(B \cup
  \{v_0 - y\})}$ (these are equal because $u_0 - v_0, x - y$ are
both in $S$).  Note that $x, y, v_0, u_0$ are all in $S'$. 
By Lemma~\ref{dimlem}, $\dim{(\mathcal{F}_{u_0} \cap S')} +
\dim{(\mathcal{F}_{v_0} \cap S')} \leq \dim{S'} - 1$. 
By the choice
of the vector $u_0 - x$ which we used to extend $S$ to $S'$, and the
fact that $x \in \mathcal{F}$, we have that $\mathcal{F}_{x} \cap S' =
\mathcal{F}_{x} \cap \mathcal{F} \cap S' = \mathcal{F}_{x} \cap S$.  Likewise
$\mathcal{F}_{y} \cap S'=\mathcal{F}_{y} \cap S$.

Since all the faces of $P$ are closed, we may assume that $u_0$, $v_0$ are close enough to $x$, $y$ so that 
$\mathcal{F}_{x} \subseteq \mathcal{F}_{u_0} \cap \mathcal{F}$ and $\mathcal{F}_{y} \subseteq
\mathcal{F}_{v_0} \cap \mathcal{F}$. 
Now $x,y \in \mathcal{F}$ but none of the points $u_n$, $v_n$ can
lie in $\mathcal{F}$, so in particular, $u_0 - x$ is not in the subspace 
defined by the affine extension of $\mathcal{F}$ (that is, the translation of the affine space which contains
$\vec 0$), and so
$\dim{(\mathcal{F}_{x} \cap S)} \leq \dim{(\mathcal{F}_{u_0} \cap \mathcal{F}
  \cap S')} < \dim{(\mathcal{F}_{u_0} \cap S')}$.
Likewise
$\dim{(\mathcal{F}_{y} \cap S)} < \dim{(\mathcal{F}_{v_0} \cap
  S')}$, and so overall we have $\dim{(\mathcal{F}_{x} \cap
  S)}+\dim{(\mathcal{F}_{y} \cap S)} \leq \dim{(\mathcal{F}_{u_0}
  \cap S')} + \dim{(\mathcal{F}_{v_0} \cap S')} -2 \leq \dim{S'} - 3
= \dim{S} - 2$.
\end{proof}

\begin{lemma} \label{dimlemd}
Suppose $S_1, S_2$ are affine linear subspaces of $\mathbb{R}^d$.  
Let $T$ be the affine extension of $S_1 \cup S_2$. Then $\dim{T} \leq \dim{S_1} + \dim{S_2} + 1$
\end{lemma}

\begin{proof}
We may assume without loss of generality that $S_1$ contains the
origin.  Let $B_1$ be a basis for $S_1$.  Let $B_2$ be a basis for the
translate of $S_2$ which contains the origin, and let $x \in S_2$.  We
claim $T' = \lspan{(B_1 \cup B_2 \cup \{x\})}$ contains $S_1 \cup
S_2$.  Indeed, any vector $v \in S_1$ is in $\lspan{B_1}$, whereas any
vector $v \in S_2$ has $v-x \in \lspan{B_2}$, and so $v \in \lspan(B_2
\cup \setof{x})$.  Since $T \subseteq T'$, and $\dim{T'} \leq
\dim{S_1} + \dim{S_2} + 1$, this gives the required estimate on the
dimension.
\end{proof}

The above analysis of the critical line segments now allows us to give the proof
of the following Lemma~\ref{dimleme}, which we use for the proof of 
the cone lemma, Lemma~\ref{conepart}. 

\begin{lemma}  \label{dimleme}
Let $P$ be a $d$-polytope in $\R^d$, and let $S\subseteq \R^d$ be an affine subspace, 
with $\dim(S)\geq 2$. Then there are finitely many co-dimension $1$ affine subspaces
$S_1,\dots,S_k \subseteq S$ such that if $x,y \in S$ and $\ell_{x,y}$ is critical for 
some intersection $\mathcal{F}$ of faces of $P$, then for some $i$ we have $y\in S_i(x)$. 
Recall $S_i(x)=S_i+ (x-v_i)$
denotes the translation of $S_i$ containing $x$.
\end{lemma}

\begin{proof}
Consider all possible pairs $(A,B)$ where $A,B$ are subsets of the
faces of $P$ and if $A', B'$ denote the corresponding sets of
hyperplanes then $\dim (\bigcap A' \cap S(\vec 0)) + \dim( \bigcap
B'\cap S(\vec 0))\leq \dim(S)-2$.  For each such pair $(A,B)$, let $x$
be any point in $\bigcap A$, and note that by Lemma~\ref{dimlemd}, if
$T$ is the smallest linear subspace containing both $(\bigcap A - x)
\cap S(\vec 0)$ and $(\bigcap B - x) \cap S(\vec 0)$, then $\dim(T)
\leq \dim(S) - 1$.  The resulting space $T(A, B) = T$ does not depend
on the choice of $x$.

By the proof of Lemma~\ref{translationb} we may assume that $\ell_{x,y}$ is extreme. 
Now by Lemma~\ref{dimlemc}, if $\ell_{x, y}$ is extreme and critical for an intersection $\mathcal{F}$
of faces of $P$, then 
if $A$ is the set of faces of $P$ which $x$ lies on and $B$
is the set of faces of $P$ that $y$ lies on, then $(A, B)$ is a pair
satisfying $\dim (\bigcap A' \cap S(\vec 0)) + \dim( \bigcap B'\cap
S(\vec 0))\leq \dim(S)-2$, and so $T(A, B) \subseteq S(\vec 0)$ has
co-dimension at least $1$ in $S(\vec 0)$. We also have that $y-x
\in T(A, B)$, and so $y \in T(A, B)(x)$, so that the collection of
finitely many $T(A, B)$, translated to lie inside of $S$ are as desired.
\end{proof}

\begin{lemma} \label{notcrit}
Let $P$ be a $d$-polytope in $\R^d$, and let $S\subseteq \R^d$ be an affine subspace. 
Suppose $x, y\in P$ and $\ell_{x,y}$ is extreme and maximal 
for $P$ with $\ell_{x,y}$ lying in exactly the set $\mathfrak{F}$ of faces of $P$. 
Suppose $\ell_{x,y}$ is not critical for $\mathcal{F}=\cap \mathfrak{F}$. 
Then there are neighborhoods $U,V$ of $x,y$ respectively such that for any $u \in 
U$, $v \in V$, if $\ell_{u,v}$ can be translated and scaled to be in $\mathcal{F}$ 
then $\ell_{u,v}$ can be translated and scaled to be extreme and maximal 
for $P$ and lie in exactly the faces $\mathfrak{F}$.
\end{lemma}

\begin{proof}
Since $\ell_{x,y}$ is not critical for $\mathcal{F}$, there are neighborhood $U$, $V$
of $x$, $y$ such that for any $u \in U$, $v\in V$, if $\ell_{u,v}$ can be scaled and
translated to be in $\mathcal{F}$, then $\ell_{u,v}$ can be scaled and translated 
to be maximal in $P$ and lying in $\mathcal{F}$. Furthermore, Lemma~\ref{translationa} 
showed that we may take $U$, $V$ so that for such $u \in U$, $v \in V$, we may 
find a maximal copy of $\ell_{u,v}\subseteq \mathcal{F}$ in small enough neighborhoods of $x$, $y$ so that 
$u$ does not lie on any faces that $x$ does not lie on, and likewise for $v$. 
Also, by Lemma~\ref{translationb} we may further assume that $\ell_{u,v}$
is extreme and maximal. Also, $\ell_{u,v}\subseteq \mathcal{F}$ (as translating a maximal segment to make
it extreme does not decrease the set of faces it lies on) but also from the 
assumption on $U,$ $V$ we have that $\ell_{u,v}$ does not lie on more faces than $\ell_{x,y}$. 
\end{proof}

Now we are ready to prove the generalized cone lemma, Lemma~\ref{conepart},
which we  restate. 
\begin{lemma*} 
Let $P$ be a polytope in $\R^d$. Suppose $A \subseteq S$ is uncountable where $S\subseteq \R^d$ is an affine subspace,
and assume that every hyperplane in $S$ intersects $A$ in a countable set.
Suppose that $H_1, \dots H_k$ are distinct hyperplanes in $S$ inducing the family of cones $C_1
\dots C_\ell$. 
Then there is a cone $C_i$, a subset $\mathfrak{F}$ of the faces of $P$  and an uncountable subset
$A' \subseteq A$ so that for every uncountable $B\subseteq A'$ there is an uncountable $D\subseteq B$
such that for every point $x \in D$, $x$ is a strong limit
both from $C_i^+(x) \cap D$ and $C_i^-(x) \cap D$ of points $y$ in $G(x,\mathfrak{F})$. 
\end{lemma*}

\begin{proof}
We first note that it suffices to prove the following weaker version of the theorem:
for any uncountable $A\subseteq S$ and cones $C_1,\dots,C_\ell$, 
there is an uncountable $A'\subseteq A$, a $C_i$ and a collection of faces $\mathfrak{F}$ such that 
every $x \in A'$ is a strong limit point of $C_i^+(x)\cap A'$ and 
and a strong limit point of $C_i^-(x)\cap A'$ of points $y$ such that there is a scaled, translated copy $P'$
of $P$ with $\ell_{x,y}$ extreme and maximal in $P'$ and for every face $F'$ of $P'$, $\ell_{x,y}\subseteq F'$ iff 
the corresponding face $F$ of $P$ is in $\mathfrak{F}$. 
To see this, suppose the conclusion of Lemma~\ref{conepart} fails.
We let $G(x,\mathfrak{F})$ be the set 
of such points $y$ as in the statement of the theorem for a given $x$ and set of faces $\mathfrak{F}$. 
Let $(C_{p_i},\mathfrak{F}_{q_i})$, for $1\leq i \leq m$,
enumerate the pairs of cones and subsets of faces. 
Let $A_0=A$,
and supposing $A_{i-1}$ has been defined, let $A_i\subseteq A_{i-1}$ be uncountable
so that for all uncountable $E\subseteq A_i$ it is not the case that 
for every $x \in E$ that $x$ is a strong limit
both from $C_{p_i}^+(x) \cap E$ and  $C_{p_i}^-(x) \cap E$ of points $y\in G(x,\mathfrak{F_{q_i}})$.
We can get $A_i$ from the assumption that the theorem fails for the pair $(C_{p_i},\mathfrak{F}_{q_i})$ and $A_{i-1}$. 
Then $A_{m}$ violates the weaker version of the theorem.

To prove the weaker version, fix the uncountable set $A\subseteq S$, the hyperplanes $H_1,\dots,H_k$
and corresponding cones $C_1,\dots,C_\ell$.  We may assume that 
the hyperplanes $H_1,\dots,H_k$ includes the hyperplanes from Lemma~\ref{dimleme}. 

As in the proof of Theorem~\ref{cone}, we define an operation thinning
the uncountable set $A$ out, but in this case we use the enumeration
of the pairs $(C_{p_i}, \mathfrak{F}_{q_i})$ to do so.  Let $A_0=A$.
Consider the first pair $(C_{p_1}, \mathfrak{F}_{q_1})$. We ask if
there is an uncountable $A'_0 \subseteq A_0$ such that for every $x
\in A'_0$ there is a neighborhood $U$ of $x$ such that either the set
$U\cap C^+_{p_1}(x)\cap A'_0\cap G(x,\mathfrak{F}_{q_1})$ 
or the set 
$U\cap C^-_{p_1}(x)\cap A'_0 \cap G(x,\mathfrak{F}_{q_1})$  is countable.  
If this is the case then we may thin
$A'_0$ out to an uncountable set $A_1 \subseteq A'_0$ and fix a symbol
$s(1)\in \{+, -\}$ so that for all $x \in A_1$, there is a
neighborhood $U$ of $x$ such that $U \cap
C^{s(1)}_{p_1}(x) \cap A_1\cap G(x,\mathfrak{F}_{q_1})$ is countable.
That is, by thinning from $A'_0$ to $A_1$ we fix the side
of the cone $C_{p_1}$ which is weakly isolating $x$ in $A_1$ with
respect this property.  We continue in this manner for $m$ stages,
assuming that we can at each stage go from $A_k$ to $A'_k\subseteq
A_k$ using the pair $(C_{p_k},\mathfrak{F}_{q_k})$. We then define $A_{k+1} \subseteq A'_k$ and
$s(k)$ exactly as in the first step. Assuming this construction
continues for all $m$ steps, we end with an uncountable set $A_m
\subseteq A$ and for each pair $(C_{p_i}, \mathfrak{F}_{q_i})$, we
have a side $s(i)$ of the cone $C_{p_i}$.

In this case, for each $x \in A_m$ we let $B_x$ be a rational ball
containing $x$ such that $\forall i$, the set $B_x \cap
C^{s(i)}_{p_i}(x) \cap A_m \cap G(x,\mathfrak{F}_{q_i})$ 
is countable. We may fix a rational ball
$B$ so that $B_x=B$ for an uncountable subset $A'_m$ of $A_m$.  Let
$x,y \in A'_m$ be distinct strong limit points of $A'_m$, so in
particular $x,y \in B=B_x=B_y$.  Since the set of points of $A'_m$
which are strong limit points of $A'_m$ is uncountable, from our
assumption on $A$ we may assume that the line $\ell_{x,y}$ between $x$
and $y$ does not lie on any of the $H_1,\dots,H_k$, and in particular,
by Lemma~\ref{dimleme}, is not a critical line segment for any
intersection $\mathcal{F}$ of faces of any copy $P'$ of $P$ in which
it is extreme and maximal. Let $C$ be the (unique) cone such that $y \in C^{\pm}(x)$
(and so also $x \in C^{\pm}(y)$). Note that $y$ is in the interior 
of $C^{\pm}(x)$, and so also $x$ is in the interior of $C^{\pm}(y)$. 
Let $\mathfrak{F}$ be a subset of the faces of $P$  
such that some extreme maximal (translated and scaled) copy of $\ell_{x,y}$
lies on exactly the faces $\mathfrak{F}$ of $P$. So, $y  \in G(x,\mathfrak{F})$
(and so also $x \in G(y,\mathfrak{F})$). 
Let $i$ be such that $(C, \mathfrak{F})=(C_{p_i},\mathfrak{F}_{q_i})$. 
Without loss of generality we may assume $y \in C^{s(i)}(x)$. 
Let $V$ be a neighborhood of $y$ contained in the interior of 
$B\cap C^{s(i)}(x)$. On the one hand, since $V\subseteq B$ there are only countably 
many $y' \in V$ such that $y' \in G(x, \mathfrak{F})$. On the other hand, 
from Lemma~\ref{dimleme}, since $\ell_{x,y}$ is extreme and maximal in some 
scaled translated copy $P'$ of $P$ with $\ell_{x,y}$ 
lying on exactly the faces of $P'$ corresponding to $\mathfrak{F}$ and since 
$\ell_{x,y}$ is not critical for $\mathcal{F}=\cap \mathfrak{F}$, by Lemma~\ref{notcrit}
we have that for all $y'$ in a small enough neighborhood
of $y$ that $y' \in G(x, \mathfrak{F})$. This contradiction completes
the proof.

\end{proof}

\section{Compact, convex sets in $\R^2$} \label{sec:r2}

In this section we consider the case where $P$ is a compact, convex (not necessarily
strictly convex) set in $\R^2$ with non-empty interior, and prove the corresponding version of our main result,
Theorem~\ref{thmc}. For the rest of this section we fix the compact, convex set $P
\subseteq \R^2$. The argument in this case will require a generalization of the cone lemma
Lemma~\ref{cone}.

Since $P$ is not necessarily strictly convex, Lemma~\ref{lemma:delta} does not apply and there may be
points $x \in \partial P$ with $\delta_P(x)=0$. The following lemma will serve as a replacement
for Lemma~\ref{lemma:delta} in the current case.

\begin{lemma} \label{lemma:r2delta}
Let $P$ be a compact, convex set set $\R^2$. Then there is a $\delta>0$ and a finite set
$F\subseteq \partial P$ such that $\delta_P(x)>\delta$ for all $x \in \partial P\sm F$, and
$\delta_P(x)=0$ for all $x \in F$. 
\end{lemma}

\begin{proof}
If the conclusion fails then there is a sequence of distinct points $\{ x_n\} \subseteq \partial P$
with $x_n \to x \in \partial P$ and with $\delta_P(x)=0$ and $\delta_P(x_n)\to 0$.
We use here the facts that $\partial P$ is compact and if $\{ x_n\} \subseteq \partial P$, $x_n \to x$,
and $\delta_P(x_n)\to 0$, then $\delta_P(x)=0$. To see this, suppose towards a contradiction that
$\{ x_n\} \subseteq \partial P$ are distinct points, $x_n \to x $ (so $x \in \partial P$),
$\delta_P(x_n) \to 0$, but $\delta_P(x)>0$. For each $n$, let $y_n \in \partial P$ be such that
$y_n \in M_P(x_n)$ and if $\ell_n$ is the vector from $x_n$ to $y_n$ then
$\frac{\ell}{\| \ell\|} \cdot H =\delta_P(x_n)$ for some supporting hyperplane (line)
for $P$ which contains $x_n$. Likewise, let $y \in \partial P$ be in $M_P(x)$.
By passing to a subsequence we may assume that $y_n \to y \in \partial P$. From
Lemma~\ref{lem:mlc} we have that $y \in M_P(x)$. From Remark~\ref{rem:lb} we have that
$x\neq y$, that is, the vector $\ell$ from $x$ to $y$ is non-degenerate. 
Since $\delta_P(x)>0$ by assumption, there is a separation between $\ell$ and any
supporting line for $P$ at $x$. This easily implies that the only points of $\ell$
which are not in $P^\circ$ are $x,y$. For we cannot have $\ell \subseteq \partial P$
and so there is a point $z \in \ell$ with $z \in P^\circ$. From the convexity of $P$ it
then follows that every point on $\ell$ except $x$ and $y$ is in $P^\circ$. 
Let $z\in \ell$ be a point distinct from $x$ and $y$. Let $B$ be a ball about $z$
with $B \subseteq P$. Then for $n$ large enough, so $x_n, y_n$ are sufficiently close
to $x,y$ we have that $\ell_n$ intersects $B$ and in fact since $\delta_P(x_n)\to 0$
we will have that for large enough $n$ that some supporting line $H_n$ for $P$
at $x_n$ will intersect $B$. This contradicts the fact that $B\subseteq P^\circ$.

So, assume $\{ x_n\}\subseteq \partial P$ are distinct points, $x_n \to x$, $\delta_P(x_n)\to 0$, and
$\delta_P(x)=0$. As above, for each $n$ let $y_n \in \partial P$ be such that
$y_n \in M_P(x_n)$ and if $\ell_n$ again denotes the vector from $x_n$ to $y_n$
then for some supporting hyperplane $H_n$ for $P$ at $x_n$ we have that
$\frac{\ell_n}{\| \ell_n\|} \cdot H_n=\delta_P(x_n)$. As before we may assume $y_n \to y$
for some $y \in \partial P$. So $y \in M_P(x)$, and 
the above argument shows that $\ell \subseteq \partial P$, where $\ell$ is the vector from
$x$ to $y$. Also, $x\neq y$ so $\ell$ is non-degenerate. Consider a sufficiently large $n$,
so $x_n,y_n$ are sufficiently close to $x,y$. Let $z_n$ be the midpoint of $\ell_n$.
We cannot have $z_n$ on the line segment $\ell$ from $x$ to $y$ as otherwise
if $x_n\notin \ell$ then we contradict that the line determined by $\ell$ is a supporting line for $P$,
and if $x_n \in \ell$, then by maximality of $\ell$ and $\ell_n$ we would have $\ell=\ell_n$
and so $x=x_n$ which could only happen for one $n$. 
Thus, the triangle $T=T(z_n,x,y)$ determined by the points $z_n$, $x$, $y$
is non-degenerate. By convexity $T\subseteq P$. For $m$ sufficiently large, $\ell_m
\cap T^\circ \neq \emptyset$. For $m$ sufficiently large, if $x_m \notin \ell$, then 
we have $H_m \cap T^\circ \neq \emptyset$, which contradicts the fact that
$H_m$ is a supporting line for $P$ and $T^\circ \subseteq P$. 

So, we may suppose that for large enough $n$ we have $x_n \in \ell$, and so $y_n \notin \ell$. 
Fix an  $n$ and consider the cone $C$ from $x$ with sides $\ell$ and $\ell_{x,y_n}$. 
For $n'>n$ large enough, we have that $y_{n'}$ lies in $\interior{C}$. Note that by convexity the convex hull 
$Q$ of the points $y,x,y_n,y_{n'}$ is contained in $P$. Since $y_{n'}\in \interior{C}$, the ray
from $y_{n'}$ parallel to $\ell$ in the direction towards $x$ contains a non-degenerate interval 
inside  $Q$. Thus the line segment $\ell_{n'}$ can be translated in the direction $x-y$ to 
stay inside $P$ and with an endpoint in $\interior{P}$. This contradicts the maximality of $\ell_{n'}$.

\end{proof}

Given the compact, convex set $P \subseteq \R^2$, let $\delta>0$ and $F\subseteq \partial P$
be from Lemma~\ref{lemma:r2delta}. Consider a point $x \in F$, so $\delta_P(x)=0$.
Let $L$ be a supporting hyperplane (line) for $P$ at $x$. We may identify $L$ with the $x$-axis
for convenience of notation, with $x$ identified with the origin and with $P$ in the upper half-space. 
Consider the left and right tangents to $P$ at $x$. For example, the right tangent
is the ray $r_\theta$ from the origin with angle $\theta$ (measured counterclockwise from the
$x$-axis as usual) where $\theta$ is largest so that $P$ lies on the half-space
which is determined the line $L_\theta$ extending $r_\theta$ and which is disjoint from
the sector determined by the positive $x$-axis and the ray $r_\theta$.
We let $r_x$ denote the unit right tangent to $P$ at $x$. 
We similarly define the unit left tangent $l_x$ to $P$ at $x$. 
We may suppose that $\ell_x \neq r_x$ as otherwise $P$ lies on a line segment in
$\R^2$, which violates the assumption that $P$ has non-empty interior.

Since $\delta_P(x)=0$ (as $x \in F$), we have that at least one of $\ell_x$, $r_x$
lies in the direction of a line segment $\ell \in M_P(x)$.

\begin{claim} \label{bnc}
Suppose that only one of these tangents, say $r_x$,
lies in the direction of a vector in $M_P(x)$.
There is an $\eta_x>0$ such that if the directed line segment $\ell$ from $x$
to a point $y \in \partial P$ lies in $M_P(x)$ and $\ell$ is not coincident with
$r_x$, then $|\frac{\ell}{\| \ell\|} \cdot N(\ell_x)| >\eta_x$,
where $N(\ell_x)$ is the normal vector to the vector $\ell_x$.
\end{claim}

\begin{proof}
If the claim fails then there is a sequence of points $\{ y_n\}\subseteq \partial P$
with the line segments $\ell_n$ from $x$ to $y_n$ in $M_P(x)$ and
such that $|\frac{\ell_n}{\| \ell_n\|} \cdot N(\ell_x)| \to 0$. The lengths of the
line segments $\ell_n$ are bounded away from $0$, and we may assume that
$y_n \to y \in \partial P$ with $x \neq y$. By convexity, the line segment
$\ell$ from $x$ to $y$ lies in $P$. We also have that
$|\frac{\ell}{\| \ell\|} \cdot N(\ell_x)|=0$, and so $\ell$ is parallel to $\ell_x$.
From the definition of $\ell_x$ we easily have that $\ell$ points in the same direction as
$\ell_x$ (as $\ell_x$, $\ell$ both point into the same half-space determined by the
line $L$ as in the definition of $\ell_x$ above). From Lemma~\ref{lem:mlc}
we have that $\ell \in M_P(x)$. This contradicts our assumption that
$\ell_x$ does not lie in the direction of a vector in $M_P(x)$.
\end{proof}

Let $F'\subseteq F$ be the points $x \in F$ such that  exactly 
one of the $\ell_x$, $r_x$ directions has a maximal line segment with endpoint $x$
which lies in $\partial P$. 
Let $\eta=\min \{ \eta_x \colon x \in F'\}$, and let $\delta'=\min \{ \delta,\eta\}$. 
Let $H_1,\dots,H_k$ be a finite set of hyperplanes (lines) in $\R^2$ containing all of the $\ell_x$ and $r_x$
(translated to pass through the origin)  for each $x \in F$, and also such that if $C_1,\cdots,C_\ell$
is the set of corresponding cones (in this case $\ell=k$) then the angle of each cone
$C_j$ is less than $\delta'$.

Let $A \subseteq \R^2$ be such that $\II$ has a winning strategy $\tau$ in the no-$\beta$ McMullen
game for $A$. Suppose towards a contradiction that $A$ is uncountable.
As in the proof of Theorem~\ref{thma} we may assume that $A\cap \ell$ is countable
for every line $\ell \subseteq \R^2$ (as otherwise we may replace $A$ with $A\cap \ell$,
and the argument is a much easier version of the argument to follow). 
Let $N$ be a countable elementary substructure (of a large $V_\kappa$)
containing $\tau$ and $H_1,\dots,H_k$ (and so the cones $C_1,\dots,C_k$). We define the trees
$T$, $T^N$, and $T_x$, $T^N_x$ as before. As before, we get a position $p=(B_0, B_1,\dots, B_{n+1})$ (so $n$ is even)
in the tree $T^N$ such that for all $x$ is some uncountable $A' \subseteq A$, $p$ is terminal in $T^N_x$. 
As before we get a rational ball $D \subseteq B_n^\circ   \sm B_{n+1}$ such that
$A'\cap D$ is uncountable and for any $x,y \in D$, $\bp{x}{y} \subseteq B_n^\circ \sm B_{n+1}$.

We apply Lemma~\ref{cone} to $A' \cap D$ and the hyperplanes $H_1,\dots,H_k$. 
Let $A'' \subseteq A'$ be the uncountable set and $C_i \in \{ C_1,\dots, C_k\}$ be the cone
from Lemma~\ref{cone}.

Let $x \in A''$. We have that $x$ is a strong limit point of $A'' \cap C^+_i(x)$ and $A'' \cap C^{-}_i(x)$,
and since $A\cap \ell$ is countable for every line $\ell$, we can choose $y \in A'' $ with $y$ in the interior
of $C^+_{i}(x)$. So, $x$ is in the interior of $C^{-}_i(y)$. Consider $\bp{x}{y}$.
By definition the line-segment $\ell_{x,y}$ from $x$ to $y$ is maximal in $\bp{x}{y}$.
Without loss of generality assume $x \notin \tau( p\concat \bp{x}{y})$. Let $B_1$ be a rational ball
containing $x$ with 
$B_1 \cap \tau(\bp{x}{y})=\emptyset$.
If $\delta_{\bp{x}{y}}(x)\neq 0$, then since $\delta' \leq \delta$ and the cone angle of $C_i$
is $\leq \delta'$ it follows that if $z$ in the interior of $C_i^+(x)$ is close enough to $x$ then $z \in \bp{x}{y}^\circ$.
Suppose now that $\delta_{\bp{x}{y}}(x)=0$. Note that $y \notin \ell_x \cup r_x$. 
We consider cases as to whether both or just one of $\ell_x$, $r_x$ is contained in $\partial P$.
In the case that both are contained in $\partial P$ 
we claim that the cone
$C^+_i(x)$ is contained in the cone with vertex at $x$ and sides given by $\ell_x$, $r_x$.
This follows from the facts that $\ell_{x,y} \subseteq C^+_i(x)$, and the lines $\ell$ and $r$
through the origin parallel to $\ell_x$ and $r_x$ were included among the $H_1,\dots,H_k$. 
On the other hand, if either $\ell_x \cap \bp{x}{y}=\{x\}$ or $r_x \cap \bp{x}{y}=\{x\}$,
then from Claim~\ref{bnc} we have that the angle between $\ell_{x,y}$ and $\ell_x$ (or $r_x$
respectively) is at least $\eta_x \geq \delta'$. Thus, in both cases if $z$ is in the interior of
$C^+_i(x)$ is close enough to $x$, then $z \in \bp{x}{y}^\circ$.
Let $z\in A''$ be such a point close enough to $x$.

So, in all cases we have  $z \in A'' \cap \bp{x}{y}^\circ$. 
Let $B_2\subseteq B_1\cap \bp{x}{y}^\circ$
be a rational ball containing $z$. The argument now finishes exactly as in Theorems~\ref{thma} and \ref{thmb}.

\section{Ordinal Analysis} \label{sec:oa}

Throughout this section we fix a compact, convex set $P\subseteq \R^d$ for which
the no-$\beta$ McMullen game is equivalent to the perfect set game. By Theorems~\ref{thma},
\ref{thmb} and \ref{thmc}, this includes all strictly convex sets in $\R^d$, all compact, convex
sets in $\R^2$, and all polytopes in $\R^d$.

Recall from Definition~\ref{derivative} of \S\ref{sec:intro} that for a 
set $A\subseteq \R^d$ we have defined a derivative
notion $A'_P \subseteq A$. Iterating this defined the sets $A^\alpha_P$.
This process stops at a countable ordinal $\alpha$, and we set $A^\infty_P=A^\alpha_P$.

Our geometric consequences result, Theorem~\ref{geomthm},  
follows immediately from Theorems~\ref{thma}, \ref{thmb}, and \ref{thmc}
and the following Theorem~\ref{thm:derwin}, which relates the ordinal analysis of the derivative
for closed sets to the no-$\beta$ McMullen game.

We note that (\ref{der1}) of the following theorem would follow easily if in Definition~\ref{derivative}
we did not have the clause, in the definition of a good copy, involving limits along hyperplanes $H$ 
but only consider limits from points in the interior of the copy of $P$. For in this case, $\I$ could win
by playing good copies for $A^\infty_P$. However, this clause of Definition~\ref{derivative}
is necessary for the proof of the following theorem in the polytope case (but not in the
strictly convex case).

\begin{theorem} \label{thm:derwin}
Let $P\subseteq \R^d$ be a compact, convex set and assume either $P$ is strictly convex, or $d=2$, or $P$ is a
polytope.
Let $A\subseteq \mathbb{R}^d$ be a closed set and consider the sets $A^\alpha_P$, and $A^\infty_P$ defined above.
\begin{enumerate} 
\item \label{der1}
If $A^\infty_P$ is nonempty, then player $\I$ has a winning strategy in the
no-$\beta$ McMullen game $G_P(A)$.
\item \label{der2}
If $A^\infty_P$ is empty, then player $\II$ has a winning strategy in the
no-$\beta$ McMullen game $G_P(A)$. 
\end{enumerate}
\end{theorem}

\begin{proof}
We first show that (\ref{der1}) follows from (\ref{der2}). Assuming (\ref{der2}), let
$A\subseteq \mathbb{R}^d$ be a closed set with $A^\infty_P \neq \emptyset$. 
This implies that $A^\infty\neq \emptyset$, where we recall that $A^\infty$ refers to the ordinary
Cantor-Bendixson derivative. Since $A$ is closed, this means that $A$ is uncountable and in fact
$A^\infty$ is a perfect set. By Theorems~\ref{thma}, \ref{thmb}, and \ref{thmc}, $\II$
does not have a winning strategy in the game $G_P(A)$. The game $G_P(A)$ is a closed game for $\I$,
and so is determined by Gale-Stewart (see \cite{GaleStewart1953} or \S 6A of \cite{MoschovakisBook}).
So, $\I$ has a winning strategy for $G_P(A)$.

Now we show that (\ref{der2}) holds. 
Assume that $A^\infty_P=\emptyset$.
Since $P$ is fixed for the rest of the argument and we have no need for the ordinary Cantor-Bendixson derivative,
we drop the subscripts $P$ for the rest of the argument,
so we write $A^\infty$ for $A^\infty_P$ and $A^\alpha$ for $A^\alpha_P$. 
If $A^\infty$ is empty, and player $\I$ plays some copy $P_0$, then
every point $x$ of $A$ on the boundary of $P_0$ has some rank $\alpha$
so that $x \in A^\alpha \setminus A^{\alpha+1}$.

Let $\alpha_0$ be minimal so that $P_0 \cap A_{\alpha_0+1} =
\emptyset$. Such an $\alpha_0$ exists: first note that if there are no
points of $A=A_0$ on the boundary of $P_0$, then we can immediately
win the game, since by a compactness argument, there are no points of
$A$ within some fixed $\epsilon$ of the boundary. Next note that each
of the $A_\alpha$ are closed, and so $A_\alpha \cap \partial P_0$ is
compact, and so the minimal $\alpha$ so that $\partial P_0 \cap
A_\alpha = \emptyset$ cannot be a limit ordinal, as we would have a
sequence of descending nonempty compact sets whose intersection is
empty.

The next claim will essentially finish the proof by allowing us to construct a 
decreasing sequence of ordinals starting from $\alpha_0$, which will give a strategy for $\II$.

\begin{claim} \label{mainclaim}
Suppose $A_{\alpha+1}\cap P_n=\emptyset$ for $P_n$
some latest ball by player $\I$ in the no-$\beta$ 
McMullen game.  There is a finite-round strategy by player $\II$ which
eliminates all of $A_\alpha$.
\end{claim}

\begin{proof}
Suppose $A_{\alpha+1}\cap P_n = \emptyset$.  Note $A_\alpha \cap
\partial P_n$ is compact, and every point in it is not in
$A_{\alpha+1}$, so there is some rational ball around each point
which contains no good copy of $P$ for $A_\alpha$. 
For each point in $\partial P_n
\setminus A_\alpha$, there is some rational neighborhood which
contains no points of $A_\alpha$.  By the compactness of $\partial P_n$,
there is some finite subcover, which covers some closed
$\epsilon$-neighborhood of $\partial P_n$.  This closed strip is
compact and so by the Lebesgue covering lemma, there's some $\delta>0$
small enough so that any ball of diameter less than $\delta$ inside this
$\epsilon$-neighborhood must be
contained in one of these rational neighborhoods,
and thus contain no good copy of $P$ for $A_\alpha$.
Play a move $P_{n+1}$ as  $\II$ to
delete enough of $P_n$ so that any move played by $\I$ must lie in this
neighborhood of $\partial P_n$ and must be of diameter less than $\delta$.
Now if player $\I$ plays $P_{n+2}$, then $P_{n+2}$
is not a good copy of $P$ for $A_\alpha$.

\begin{lemma} \label{clem}
Let $P$ be a compact, convex set in $\R^d$ with non-empty interior.
Let $B\subseteq \partial P$ be compact and such that
there is a smaller copy of $P$, that is $Q=s+tP$ for $t<1$, with $B\subseteq Q$.
Then there is a neighborhood $U$ of $B$  so that for any $\epsilon \in (0,1)$,
$U\cap P$ is in the interior of $\epsilon s+P$.
\end{lemma}

\begin{proof}
Note that $-s+B \subseteq -s+Q = tP \subseteq \interior{P}$. Let $V$ be a neighborhood of
$-s+B$ with $-s+B \subseteq V\subseteq P$. Let $U=s+V$, so $B\subseteq U$.
Let $\epsilon \in (0,1)$. Let $x \in U\cap P$. Consider the line segment $L$ between
$x$ and $-s+x \in V \subseteq P$. By convexity, $L\subseteq P$. Also by convexity, since
$-s+x \in \interior{P}$ we have that all points of $L$ except possibly $x$ are in $\interior{P}$.
In particular $-\epsilon s +x \in \interior{P}$. Thus $x$ is in the interior of $\epsilon s +P$.
So, $U\cap P$ is in the interior of $\epsilon s+P$ for any $\epsilon \in (0,1)$.  

\end{proof}

Let $B \subseteq A_\alpha \cap \partial P_{n+2}$
the set of points which are limit points of either $A_\alpha \cap \interior{P}_{n+2}$
or $A_\alpha \cap \partial P_{n+2} \cap H$ for some supporting hyperplane $H$ in $\R^d$.
Since $P_{n+2}$ is not good for $A_\alpha$ there is a legal move $Q$
for $\II$ which covers $B$. So, $Q= s+tP_{n+2}$ where $t<1$.
Let $U$ be a neighborhood of $B$ as in Lemma~\ref{clem} for $P_{n+2}$.
Let $V$ be a neighborhood of $B$ with $B\subseteq V \subseteq \closure{V} \subseteq U$. 
Let $C= A_\alpha \cap \interior{P}_{n+2}\sm V$. By a simple compactness argument, there is an
$\epsilon>0$ such that $C\cap N(\partial P_{n+2},\epsilon)=\emptyset$, where
$N(\partial P_{n+2},\epsilon)=\{ x \colon \rho(x, \partial P_{n+2})<\epsilon\}$ denotes the 
$\epsilon$-neighborhood of $\partial P_{n+2}$. So $C$ is actually closed. 
It follows that there is a $\delta >0$ such that if we shrink and translate $P_{n+2}$
by less than $\delta$, then that copy of $P_{n+2}$ will still contain $C$.
By Lemma~\ref{clem} and the definition of $\delta$, there is a $0 <\delta'<\delta$
and a $\delta'$-translation $P'$ of $P_{n+2}$ such that $P'$ contains
$U\cap P$ and hence $\closure{V} \cap P$. By compactness, there is a $\delta'' <\delta'$
such that the $\delta''$ shrinking of $P'$, say $P''$, also contains
$\closure{V}\cap P$. So, $P''$ contains all of $B$ as well as $A_\alpha \cap
\interior{P}_{n+2}$. It follows for the definition of $B$ that $A_\alpha \cap (P_{n+2}\sm P'')$
is a subset of $\partial P_{n+2}$ and consists of points which are 
not limit points of $A_\alpha \cap \partial P_{n+2} \cap H$ for any hyperplane $H$ in $\R^d$.
Have $\II$ play the move $P_{n+3}=P''$, and let $\I$ respond by $P_{n+4}
\subseteq P_{n+2}\sm P_{n+3}$. We have $A_\alpha \cap P_{n+4}\subseteq \partial P_{n+2}
\cap \partial P_{n+4}$.

\begin{lemma} \label{dlem}
Let $P$ be a compact, convex set in $\R^d$ with non-empty interior.
Let $Q=s+tP$ where $t<1$ be a smaller, translated copy of $P$, and assume
$Q\subseteq P$ and $\partial Q\cap \partial P\neq \emptyset$.
Then the following hold. If $P$ is strictly convex, then
$| \partial Q\cap \partial P|=1$. If $d=2$, then
$\partial Q\cap \partial P$ is contained in the union of two line segments
which are supporting hyperplanes for $P$.
\end{lemma}

\begin{proof}
Suppose $Q\subseteq P$ and $x \in \partial P\cap \partial Q$.
Let $x' =s+tx$ be the point ``corresponding'' to $x$ on $\partial Q$. 
Let $H(x,P)$ denote the set of hyperplanes for $P$ at $x$, that is,
$x+H(x,P)$ is the set of (affine) supporting hyperplanes for $P$ at $x$.
We clearly have that $H(x,P)=H(x',Q)$. Let $H\in H(x,P)$. Then $x+H$
is a supporting hyperplane for $P$ at $x$ and $x'+H$ is a supporting hyperplane
for $Q$ at $x'$. We also have that the side of the $x+H$ that $P$ lies
on must be the same side of $x'+H$ that $Q$ lies on. So we must have that
$x+H=x'+H$ as $x'$ must lie on the ``good'' side of $x+H$ (as $ x' \in Q\subseteq P$)
and $x$ is in the ``good'' side of $x'+H$ (as $x \in Q$). Thus, $x' \in x+H$.
Thus $x+H=x'+H$ is  a supporting hyperplane for both $P$ and $Q$ at both $x$ and $x'$.
It follows that if $x' \neq x$, then $x+H$
is a supporting hyperplane for any point on the line segment $L_{x,x'}$
between $x$ and $x'$ for both $P$ and $Q$, and so $L_{x,x'} \subseteq \partial P \cap \partial Q$.

Consider the map $f \colon \R^d \to \R^d$ be given by
$f(x)=s+tx$. This is a contraction map, so there is a unique fixed
point $x_\infty$ for $f$. Suppose $x \in \partial P\cap \partial Q$ is any point
in the common boundary of $P$ and $Q$. Let $x_0=x$ and $x_{i+1}=f(x_i)$.
By the above we have that each line segment $L_{x_i,x_{i+1}}\subseteq
\partial P \cap \partial Q$. Also, we have that all of the line segments
$L_{x_i,x_{i+1}}$ lie on a single line $L_x$. To see this, simply note that
for any point $y$, the points $y$, $f(y)=s+ty$ and $f(f(y))=s+st+t^2y$
line on a line. So we have that the line segment $L_{x,x_\infty}$ is
contained in $\partial P\cap \partial Q$.

Suppose first that $P$ is strictly convex, and suppose $x \in \partial P \cap \partial Q$.
Then we must have $x=x_\infty$, as otherwise we have a proper line segment
$L_{x,x_\infty} \subseteq \partial P$. So, $\partial P\cap \partial Q=\{x\}$.

Suppose now $d=2$, so $P\subseteq \R^2$. Suppose there is a point $x \neq x_\infty$
with $x \in \partial P \cap \partial Q$. So, $L_{x,x_{\infty}}\subseteq \partial P \cap \partial Q$.
Suppose there is another point $y \in \partial P \cap \partial Q$ with $y$ not on the line
$\ell_x$ containing $x$ and $x_\infty$. Then we likewise have that $L_{y,x_\infty}\subseteq
\partial P \cap \partial Q$. Let $\ell_y$ be the line containing $y$ and $x_\infty$.
We then claim that $\partial P \cap \partial Q \subseteq \ell_x \cup \ell_y$.
For suppose $z \in \partial P \cap \partial Q$ with $z \notin \ell_x \cup \ell_y$.
Consider the region $R$ of points between the two rays $L_{x,x_\infty}$ and
$L_{y,x_\infty}$ emanating from $x_\infty$. Since $\ell_x$ and $\ell_y$ are supporting hyperplanes for $P$,
we have $P \subseteq R$. Since $z \notin \ell_x \cup \ell_y$, $z$ must lie in the
interior of $R$, and in fact $L_{z,x_\infty}\sm \{ x_\infty\} \subseteq \interior{R}$.
But, points on $L_{z,x_\infty}\sm \{ x_\infty\}$ close enough to $x_\infty$ are
then in the interior of $P$, a contradiction.

\end{proof}

We claim that $C=A_\alpha \cap \partial P_{n+2}\cap \partial P_{n+4}$ is finite.
In the case where $P$ is strictly convex, $|A_\alpha \cap \partial P_{n+2}\cap \partial P_{n+4}|\leq
1$ by Lemma~\ref{dlem}. Suppose $d=2$. By Lemma~\ref{dlem}, there are at most two line
segments $L_1, L_2 $ which are supporting hyperplanes for $P_{n+2}$
such that $\partial P_{n+2}\cap \partial P_{n+4} \subseteq L_1\cup L_2$. 
If $C$ were infinite, then $C\cap L_1$ or $C\cap L_2$ would be infinite and
so have a limit point, which would be a point in $B$.  Here $B$ (as before)
denotes the set of points 
which are limit points of either $A_\alpha \cap \interior{P}_{n+2}$
or $A_\alpha \cap \partial P_{n+2} \cap H$ for some supporting hyperplane $H$ in $\R^d$.
This contradicts the fact that $P_{n+4}\cap B=\emptyset$. If the case where $P$ is a polytope,
then $C$ is a subset of finitely many supporting hyperplanes for $P_{n+2}$. By the same argument,
if $C$ were infinite we would get a point of $C$ in $B$, a contradiction.

In finitely many moves player $\II$ can delete all of $C$. This results in a play
with final position $P_{n+4+2k}$ with  $P_{n+4+2k} \cap A_\alpha=\emptyset$.

This completes the proof of Claim~\ref{mainclaim}.
\end{proof}

We have now defined a strategy $\tau$ for player $\II$ such for
any run $(P_0,P_1,P_2,P_3,\dots)$ according to $\tau$ then if we let
$\alpha_{2n}$ be the largest ordinal such that $P_{2n}\cap A_{\alpha_{2n}}\neq
\emptyset$ (assuming there is an ordinal such that $P_{2n}\cap A_\alpha \neq \emptyset$,
the largest such  ordinal is well-defined by a simple compactness argument),
then we have $\alpha_0\geq \alpha_{2}\geq \cdots$, and for each $n$
such that $\alpha_{2n}>0$ there is an $m>n$ such that $\alpha_{2m}<\alpha_{2n}$.
So, following $\tau$ results in a move $P_{2n}$ with $P_{2n}\cap A=\emptyset$
(assuming $\I$ has followed the rules). Thus, $\tau$ is a winning strategy for $\II$.

This completes the proof of Theorem~\ref{thm:derwin}.

\end{proof}

\begin{remark}
This proof of Theorem~\ref{thm:derwin} does not require the axiom of choice, $\ac$. 
In the proof of (\ref{der1}) we invoked the Gale-Stewart theorem for the closed real game $G_P(A)$.
It is not hard to see that the game $G_P(A)$ is equivalent to the version $G'_P(A)$
in which $\II$ must play scaled copies $s+tP$ with $s \in \Q^d$, $t\in \Q$ (equivalent here means that
if one of the players has a winning strategy in one of the games, then that same player has a winning
strategy in the other game). Without $\ac$, Gale-Stewart gives that one of the players has a
winning {\em quasistrategy} (the interested reader can consult \cite{MoschovakisBook}
for the definition) for $G'_P(A)$. If $\II$ has a winning quasistrategy for $G'_P(A)$,
then (since $\II$'s moves are coming from a countable set) $\II$ actually has a winning strategy
for $G'_P(A)$ and thus a winning strategy for $G_P(A)$. So we get in this case that $A$ is countable,
a contradiction. So, $\I$ has a winning quasistrategy in $G'_P(A)$. But again (using $\dc$ this time)
this gives a winning strategy for $\I$ in $G_P(A)$.
\end{remark}

\section{Some examples and Questions} \label{sec:example}

In this section we present an example of a compact, convex set $P \subseteq \R^3$
for which the methods of the previous theorems do not seem to apply. 
We also present an example due to David Simmons which shows that even in the case of
$\ell_2(\R)$ with $P$ being the closed unit ball, the conclusion of Theorem~\ref{geomthm} does not hold.
These examples naturally motivate some questions which we pose.

Recall that Theorem~\ref{thmb}
established the main result for the no-$\beta$ McMullen game for the case of polytopes $P$
in $\R^3$ in particular. The proof of this relied on the fact that there were only finitely
many faces of $P$, since there these were taken as some of the hyperplanes in
Lemma~\ref{cone}. Thus, if there are infinitely many pairwise non-coplanar
maximal line segments in $\partial P$ we cannot use this argument. The authors
have shown that in the case of $P$ being a cone in $\R^3$ that the main theorem
(the analog of Theorems~\ref{thma}, \ref{thmb}, \ref{thmc})
still holds. In this case, the set of maximal line segments in $\partial P$
does not lie in the union of finitely many hyperplanes, however these line
segments are pairwise coplanar in this case. When there are infinitely many
pairwise non-coplanar line segments in $\partial P$ then our arguments seem to break down,
though we do not know of any counterexamples to our main result.
We now present an example of a compact, convex set $P\subseteq \R^3$
such that the set of  maximal line segments in $\partial P$ contains an infinite,
pairwise non-coplanar set.

Let $C$ be a cone in $\R^3$ with vertex $v$, and circular base $B$. Let $P_0$
be the part of the cone between the base $B$ and a parallel plane $H$. Let $L_0,
L_1,\dots$ be line segments with one endpoint on the circular edge of $B$, and the other
endpoint on $H \cap \partial P_0$. We choose these line segments so that 
the line $\ell_i$ extending $L_i$ passes through $v$. We choose them so
that the $L_i$ have a unique limit segment $L_\infty$ and the angles between
$\ell_i$ and $\ell_\infty$ decreases monotonically to $0$. 
There are easily
neighborhoods $U_0, U_1,\dots$ of the line segments so that for each $i$,
$U_i$ is disjoint from the convex hull of $\bigcup_{j\neq i} U_j$. In fact,
by choosing the $U_i$ sufficiently small we may assume that for each $i$ there is
a hyperplane $H_i$ containing $L_i$ such that the convex hull of
$\bigcup_{j \neq i} U_j$ lies strictly on one side of $H_i$.

We now choose line segments $L'_i \in U_i$ such that
no distinct pair $L'_i$, $L'_j$ are coplanar, and we may in fact choose them
so that their endpoints are on the circles defined by $B$ and $H$.
We may also choose the $L'_i$ close enough to the $L_i$ so that
there is a hyperplane $H'_i$ containing $L'_i$ such that the 
convex hull of
$\bigcup_{j \neq i} U_j$ lies strictly on one side of $H'_i$.
In particular, the convex hull of $\bigcup_{j \neq i} L'_j$
lies strictly on one side of $H'_i$. 
Let $P$ be the convex hull of $\bigcup_i L'_i$. We have that each of the line segments
$L'_1$ lies in the boundary of $P$, since $H'_i$ is a supporting hyperplane for $P$
containing $L'_i$.

\begin{figure}[h]
  \includegraphics[width=8cm, height=8cm]{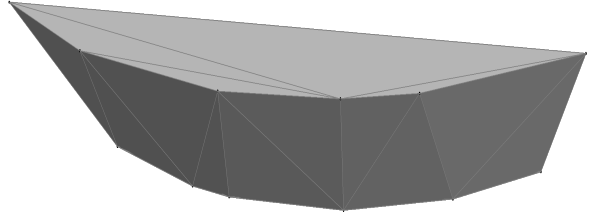}
  \vspace{-70pt}
  \caption{The convex set $P\subseteq \R^3$} \label{fig:convexhull}
\end{figure}

Note also that each of the line segments $L'_i$ is also maximal with respect to $P$
(as in Definition~\ref{def:maximal}). 
See Figure~\ref{fig:convexhull} for an illustration of the region.

There are two main results of this paper. The first is that for a compact, convex
set $P\subseteq \R^d$ satisfying the hypotheses of Theorem~\ref{thma}, \ref{thmb}, or \ref{thmc}
(that is either $P$ is strictly convex, $d=2$, or $P$ is a polytope),
then the no-$\beta$ McMullen game is equivalent to the perfect set game. The second result is that, under these
same hypotheses on $P$, for any closed set $A \subseteq \R^d$ we have that $A$
is countable iff the limiting derivative $A^\infty_P=\emptyset$.
The proofs of both of these results needed the extra hypothesis on $P$. We do not know
if either of these results is true for general compact, convex set $P\subseteq \R^d$ (for $d>2$).

\begin{question}
If $P$ is a general compact, convex set in $\R^d$, then is the no-$\beta$ McMullen game equivalent
to the perfect set game?
\end{question}

\begin{question}
If $P$ is a general compact, convex set in $\R^d$, and $A\subseteq \R^d$ a closed set, then 
is it the case that $A$ is countable iff $A^\infty_P=\emptyset$?
\end{question}

Another question one could ask is whether the results of this paper extend to infinite dimensional spaces.
We now present a result due to David Simmons which shows that in the case of
$\ell_2(\R)$ with $P$ being the closed unit ball, the results do not. In particular, 
the conclusion of Theorem~\ref{geomthm} does not hold in this case. 
The closed unit ball, of course, is not compact, so one could still ask
if the main theorem holds for compact, strictly convex sets $P$ in $\ell_2(\R)$.

\begin{example}[David Simmons]
Let $X = \setof{-1, 1}^\omega$, and let 
\[\setof{e_u \suchthat u \in \setof{-1, 1}^{<\omega}}\]
be a basis for an infinite dimensional Hilbert space $\mathcal{H}$.
Consider the map $\pi\colon X \to \mathcal{H}$
\[\pi(x)= \sum_{n \in \omega} 2^{-n} x_n e_{x \res n}\]
Let $K=\pi(X)$ and note that $K$ is compact.

\begin{claim} \label{sc}
If $y$ is in the convex hull of $K \cap \partial B(y, \rho)$ for some $\rho$, then $K \cap \interior{B(y, \rho)}=\emptyset$.
\end{claim}
\begin{proof}
Since $y$ is in the convex hull of $K \cap \partial B(y, \rho)$, there is a measure $\mu$ on $X$,
the support of which is a subset of $\pi^{-1}(K \cap \partial B(y, \rho))$, so that
\[y=\int \pi(x)~d\mu(x)\]
Note that
\begin{equation} \label{u1}
y_u=y\cdot e_u=2^{-\abs{u}-1}\left(\mu\left([u \concat 1]\right) -\mu\left([u \concat (-1)]\right)\right)
\end{equation}

Suppose for the sake of a contradiction that there was some $z \in K
\cap\interior{B(y, \rho)}$, and let $\alpha \in X$ so that $\pi(\alpha)=z$.
Since $z$ is not on the boundary of the ball, which contains the
support of $\mu$, there must be some shortest initial segment $\alpha \res (m+1)$ of $\alpha$ so
that $\mu([\alpha \res (m+1)]) = 0$.  Define $\beta$ inductively to be the
extension of $\alpha \res m$ where every initial segment has the larger
$\mu$-measure, i.e.
\[ \beta \res m = \alpha \res m\] 
and 
\[ \beta (m+k)= \begin{cases}
1 & \mu\left([ \beta \res (m-1+k) \concat 1]\right) \geq \mu\left([ \beta \res (m-1+k) \concat (-1)]\right)\\
-1 & \text{otherwise}
\end{cases}\]

Now we compute 
\[d(y, \pi(\alpha))^2 - d(y, \pi(\beta))^2\]
and we will show this quantity is nonnegative, contradicting the fact
that $\pi(\beta)$ is on the boundary of the ball, while $\pi(\alpha)$
is in the interior.  Note that the terms coming from the initial
segments $\beta \res m$ and $\alpha \res m$ cancel, and by our choice
of basis for $\mathcal{H}$, the rest of the terms are coming from
orthogonal basis vectors. Note that $y_{\alpha \res k}=0$ for all $k >m$ as $\mu([ \alpha\res k])=0$,
using Equation~(\ref{u1}).

\begingroup
\allowdisplaybreaks

\begin{align*}
&\ d(y, \pi(\alpha))^2 - d(y, \pi(\beta))^2\nonumber\\ 
=& \sum_{u} (y_u - \pi(\alpha)_u)^2 - \sum_u (y_u - \pi(\beta)_u)^2\nonumber\\
=&
[(y_{\alpha\res m}-\pi(\alpha)_{\alpha \res m})^2+ \sum_{\substack {u \leq \alpha \\ |u|>m}}
(y_u-\pi(\alpha)_u)^2 +\sum_{\substack{u \leq \beta\\ |u|>m}}
(y_u-\pi(\alpha)_u)^2 ] \\& \qquad
-[(y_{\beta\res m}-\pi(\beta)_{\beta \res m})^2+ \sum_{\substack{u \leq \alpha\\ |u|>m}}
(y_u-\pi(\beta)_u)^2 +\sum_{\substack{u \leq \beta\\ |u|>m}}
(y_u-\pi(\beta)_u)^2 ]
\\
= &(-2\alpha(m) y_{\alpha \res m}+2 \beta(m) y_{\beta\res m})
+\sum_{\substack {u \leq \alpha \\ |u|>m}} (y_u-\pi(\alpha)_u)^2
+\sum_{\substack{u \leq \beta\\ |u|>m}} (y_u)^2
\\ & \qquad
-\sum_{\substack{u \leq \alpha\\ |u|>m}} 0
-\sum_{\substack{u \leq \beta\\ |u|>m}} (y_u-\pi(\beta)_u)^2 
\\ 
\geq &
4\beta(m) y_{\beta \res m} 
+\sum_{\substack {u \leq \alpha \\ |u|>m}} (y_u-\pi(\alpha)_u)^2
-\sum_{\substack{u \leq \beta\\ |u|>m}} (y_u-\pi(\beta)_u)^2 
\\
=& 4 \beta(m)^2\mu ([\beta\res (m+1)])
+\sum_{\substack {u \leq \alpha \\ |u|>m}} (y_u-\pi(\alpha)_u)^2
-\sum_{\substack{u \leq \beta\\ |u|>m}} (y_u-\pi(\beta)_u)^2 
\\ 
\geq 
&
\sum_{n >m} (y_{\alpha \res n} - \pi(\alpha)_{\alpha \res n})^2 - \sum_{n >m}
(y_{\beta \res n} - \pi(\beta)_{\beta \res n})^2\nonumber
\\
=& \sum_{n > m} (\pi(\alpha)_{\alpha \res n})^2 - \sum_{n > m} \left(2^{-n-1}(\mu([\beta \res n \concat 1]) -
\mu([\beta \res n \concat (-1)])) - 2^{-n} \beta_n\right)^2\nonumber\\
=& \sum_{n > m} 2^{-2n} - \sum_{n > m} \left(2^{-n-1}(\mu([\beta \res n \concat 1]) -
\mu([\beta \res n \concat (-1)])) - 2^{-n} \beta_n\right)^2\nonumber\\
\geq &\ 0\ \text{since $\sgn(\beta_n) = \sgn\left(\mu([\beta \res n \concat 1]) - \mu([\beta \res n \concat (-1)])\right)$
by choice of $\beta_n$}.
\end{align*}
\endgroup
\end{proof}
\end{example}

From Claim~\ref{sc} it follows that there are no good balls in $\ell_2(\R)$ for the set $A$ (recall 
the definition of a good copy of $P$ in Definition~\ref{derivative}; here $P$ is the closed unit ball
in $\ell_2(\R)$).  
For suppose $B(y,\rho)$ were a good ball for the set $A$. Let $B\subseteq A\cap \partial B(y,\rho)$
be the points of $A\cap \partial B(y,\rho)$ which are limits of $A\cap B(y,\rho)^\circ$. 
Note here that $B(y,\rho)$ is strictly convex so no points $x$ of $A\cap \partial B(y,\rho)$
are limits of points in supporting hyperplanes for $B(y,\rho)$ containing $x$. 
Since $B(y,\rho)$ is a good ball, there is no translated smaller copy of $B(y,\rho)$ 
which contains $B$. If $y$ is in the closed convex hull $C$ of $B$, then from 
Claim~\ref{sc} we have that $A\cap B(y,\rho)^\circ=\emptyset$, and so $B$ is empty, a contradiction. 
So, $y$ is not in the closed convex hull of $B$, and thus there is a hyperplane $H$ in 
$\ell_2(\R)$ which strictly separates $C$ from an $\epsilon$ neighborhood $B(y,\epsilon)$ of $y$. 
An easy argument now shows that if $\vec n$ denotes the normal vector to $H$ 
in the direction of $C$, then for all small enough $\delta>0$, the translated ball 
$B(y+\delta \vec n, \rho)$ contains a neighborhood of $C$. So for small enough $\eta$,
$B(y+\delta \vec n, (1-\eta)\rho)$ contains $C$. This shows that $B(y,\rho)$ is not good.

So in fact we have that in $\ell_2(\R)$ there is a closed, strictly convex set $P$
(the closed unit ball) and a perfect set $A$ (the set $A=\pi(X)$ of Claim~\ref{sc})
such that $A_P'=\emptyset$. Since $P$ is not compact, this doesn't show that
Theorem~\ref{geomthm} cannot be extended to infinite dimensional spaces, but does
that the current methods encounter difficulties in trying to do so. Even supposing
the strictly convex set $P$ is compact, although some of the arguments of
Theorem~\ref{thma} go through, Lemma~\ref{cone} would no longer suffice as we would need
a version of Lemma~\ref{cone} for infinitely many cones $C_i$. 

\begin{question}
Do Theorems~\ref{thma} and \ref{geomthm} hold for strictly convex, compact $P
\subseteq \ell_2(\R)$?
\end{question}

\bibliographystyle{amsplain}
\bibliography{dst}

\providecommand{\bysame}{\leavevmode\hbox to3em{\hrulefill}\thinspace}
\providecommand{\MR}{\relax\ifhmode\unskip\space\fi MR }
\providecommand{\MRhref}[2]{%
  \href{http://www.ams.org/mathscinet-getitem?mr=#1}{#2}
}
\providecommand{\href}[2]{#2}
\begin{thebibliography}{1}

\bibitem{Davis1964}
Morton Davis, \emph{Infinite games with perfect information}, Advances in game
  theory, Annals of Math. Studies, no.~52, Princeton University Press, 1964,
  pp.~85--101.

\bibitem{FSR2016}
Lior Fishman, David Simmons, and Vanessa Reams, \emph{The
  {B}anach-{M}azur-{S}chmidt and {B}anach-{M}azur-{M}cmullen games}, Journal of
  Number Theory \textbf{167} (2016), 169--179.

\bibitem{GaleStewart1953}
David Gale and Frank~M. Stewart, \emph{Infinite games with perfect
  information}, Contributions to the theory of games, Annals of Mathematics
  Studies, no.~28, Princeton University Press, Princeton, N. J., 1953,
  pp.~245--266.

\bibitem{Martin1975}
Donald~A. Martin, \emph{Borel determinacy}, Annals of Mathematics \textbf{102}
  (1975), no.~2, 363--371.

\bibitem{Martin1985}
\bysame, \emph{A purely inductive proof of {B}orel determinacy}, Recursion
  theory (Ithaca, N.Y., 1982), Proc. Sympos. Pure Math., vol.~42, Amer. Math.
  Soc., Providence, RI, 1985, pp.~303--308.

\bibitem{McMullen2010}
Curtis~T. McMullen, \emph{Winning sets, quasiconformal maps and {D}iophantine
  approximation}, Geom. Funct. Anal. \textbf{20} (2010), no. 3, 726--740.

\bibitem{MoschovakisBook}
Yiannis~N. Moschovakis, \emph{Descriptive set theory}, Mathematical Surveys and
  Monographs, American Mathematical Society; 2nd edition, 2009.

\end{thebibliography}

\end{document}